\input amstex.tex
\documentstyle{amsppt}
\loadbold

\input epsf
\input xy
\xyoption{all}

\def\pxy{P^{1,2}_{x,y}}
\def\cixy{C^\infty_{x,y}(\bar\R,M)}
\def\czxy{C^0_{x,y}(\bar\R,M)}
\def\hgtm{H^{1,2}(\gamma^*TM)}
\def\hgu{H^{1,2}(\gamma^*V)}
\def\exy{{\Cal E}^2_{x,y}}
\def\lstm{L^2(\sigma^*TM)}
\def\lgtm{L^2(\gamma^*TM)}

\def\es{\end_{\text{sym}}}
\def\end{\operatorname{End}}
\def\supp{\operatorname{supp}}
\def\maps{\operatorname{Maps}}
\def\Cr{\operatorname{Cr}}

\def\ie{i.e\. }
\def\resp{resp\. }

\def\R{{\Bbb R}}
\def\Z{{\Bbb Z}}
\def\N{{\Bbb N}}
\def\C{{\Bbb C}}

\def\Img{\operatorname{img}}

\def\ev{\operatorname{ev}}

\def\dt{\partial_t}

\def\ddt{\tfrac\partial{\partial t}}

\def\rank{\operatorname{rank}}
\def\re{\operatorname{\frak{Re}}}
\def\Zeros{\operatorname{Zeros}}
\def\grad{\operatorname{grad}}
\def\index{\operatorname{ind}}
\def\id{\operatorname{id}}
\def\pr{\operatorname{pr}}
\def\Int{\operatorname{Int}}
\def\Pt{\operatorname{Pt}}
\def\spect{\operatorname{spect}}
\def\sm{{\text{\rm sm}}}
\def\tsm{{t,\sm}}
\def\tla{{t,\la}}
\def\la{{\text{\rm la}}}
\def\Vol{\operatorname{vol}}
\def\llangle{\langle\!\langle}
\def\rrangle{\rangle\!\rangle}

\topmatter
\title
On the topology and analysis of a closed one form. I
\\
(Novikov's theory revisited)
\endtitle

\rightheadtext{On the topology and analysis of a closed one form. I}

\thanks
Supported in part by NSF
\endthanks

\author
D.~Burghelea and S.~Haller
\endauthor
\address
The Ohio State University, Columbus, OH 43210
\endaddress
\email
burghele\@mps.ohio-state.edu
\endemail
\address
The Ohio State University, Columbus, OH 43210
\endaddress
\email
shaller\@math.ohio-state.edu
\endemail

\abstract
We consider systems $(M,\omega,g)$ with $M$ a closed smooth manifold,
$\omega$ a real valued closed one form and $g$ a Riemannian metric, so
that $(\omega,g)$ is a Morse-Smale pair, Definition~2.
We introduce a numerical
invariant $\rho(\omega,g)\in[0,\infty]$ and
improve Morse-Novikov theory by showing that the Novikov complex
comes from a cochain complex of free modules over a subring
$\Lambda'_{[\omega],\rho}$ of the Novikov ring $\Lambda_{[\omega]}$ which
admits surjective ring homomorphisms $\ev_s:\Lambda'_{[\omega],\rho}\to\C$
for any complex number $s$ whose real part is larger than $\rho$.
We extend Witten-Helffer-Sj\"ostrand results from a pair $(h,g)$
where $h$ is a Morse function to a pair $(\omega,g)$ where $\omega$ is
a Morse one form. As a consequence we show that if $\rho<\infty$ the
Novikov complex can be
entirely recovered from the spectral geometry of $(M,\omega,g)$.
\endabstract

\toc
\widestnumber\head{9}
\head 0. Introduction\endhead
\head 1. The results\endhead
\subhead 1.1 Novikov rings\endsubhead
\subhead 1.2 Morse-Smale condition and the invariant $\rho$\endsubhead
\subhead 1.3 Compactification\endsubhead
\subhead 1.4 Novikov complexes\endsubhead
\subhead 1.5 Witten-Helffer-Sj\"ostrand theory\endsubhead
\head 2. The proof of Propositions 2\endhead
\head 3. The proof of Proposition 3\endhead
\head 4. The proof of Theorem 1\endhead
\subhead 4.1 Some notations\endsubhead
\subhead 4.2 The proof for degree of rationality 1\endsubhead
\subhead 4.3 The transversality of $p$ and $s$\endsubhead
\subhead 4.4 The general case\endsubhead
\head 5. The proof of Theorem 2\endhead
\head 6. Sketch of the proof of Theorems 3 and 4\endhead
\endtoc
\endtopmatter

\document
\head 0. Introduction\endhead

Let $(M,\omega,g)$ be a system consisting of a closed connected
smooth $n$-dimensional manifold $M$, a closed one form
$\omega$ and a Riemannian metric $g$.
The form $\omega$ induces the homomorphism
$[\omega]:H_1(M;\Z)\to\R$.
Denote by $\Gamma:=H_1(M;\Z)/\ker([\omega])$.

For any two points $x,y\in M$ denote by $\Cal P(x,y)$ the set of
$\Gamma$-equivalence
classes of smooth paths $\alpha:[0,1]\to M$ with
$\alpha(0)=x$ and $\alpha(1)=y$, where we say
that $\alpha$ is $\Gamma$-equivalent to $\beta$ iff
$[\omega](\alpha\beta^{-1})=0$. Here $\alpha\beta^{-1}$ represents
the cycle obtained by going along $\alpha$ and returning along $\beta$.
The equivalence class of $\alpha$ will be denoted by $\hat\alpha$.
The juxtaposition of paths $\alpha$ and $\beta$ with $\alpha(1)=\beta(0)$
defines an action $\Gamma\times\Cal P(x,y)\to\Cal P(x,y)$ which is free and
transitive, and the obvious map $\Cal P(x,y)\times\Cal P(y,z)\to
\Cal P(x,z)$. The form $\omega$ also
associates the function $[\omega]:\Cal P(x,y)\to\R$ defined by
$[\omega](\hat\alpha):=\int_\alpha\omega\in\R$.

Suppose $\omega$ is a Morse form. Then each critical point
$x\in\Cr(\omega):=\Zeros(\omega)$ is non-degenerated and has an index,
$\index(x)\in\{0,1,\dotsc,\dim(M)\}$. The unstable set $W^-_x$ of the
vector field $X=-\grad_g\omega$, at the critical point $x$,
is the image in $M$ by the one to one immersion
$i_x^-:\R^{\index(x)}\to M$ defined in an obvious way using the trajectories
departing from $x$.
For any $x\in\Cr(\omega)$ choose an orientation $\Cal O_x$ on $\R^{\index(x)}$.
Denote by
$o:=\{\Cal O_x\mid x\in\Cr(\omega)\}$ the collection
of all these orientations.

Suppose that $(\omega,g)$ satisfies the Morse-Smale condition, cf
Definition~2 in section~1.2. For any $x\in\Cr_q(\omega)$,
$y\in\Cr_{q-1}(\omega)$ and $\hat\alpha\in\Cal P(x,y)$, S.~P.~Novikov
has associated the
integer number $I_q(x,y,\hat{\alpha})$, cf section~1.4
for definition, and has noticed the following properties:
\roster
\item
For any real number $R$ the set
$$
\bigl\{\hat\alpha\in\Cal P(x,y)\bigm| I_q(x,y,\hat\alpha)\neq 0,
[\omega](\hat\alpha)\geq R\bigr\}
$$
is finite.
\item
For any $x\in\Cr_q(\omega)$, $z\in\Cr_{q-2}(\omega)$ and
$\hat{\gamma}\in\Cal P(x,z)$
the sum
$$
\sum_{
\Sb y\in\Cr_{q-1}(\omega)\\
\hat\alpha\in\Cal P(x,y),\ \hat\beta\in\Cal P(y,z)\\
\hat\alpha\hat\beta=\hat\gamma\endSb}
I_q(x,y,\hat\alpha)I_{q-1}(y,z,\hat\beta)=0,
$$
\endroster
which means that in the sum above, the left side contains only finitely
many nonzero terms whose sum is zero.

As a consequence the collections of numbers $I_q(x,y,\hat\alpha)$ can be
algebraically organized to provide a cochain complex
$(N\bold C^*,\boldsymbol\partial^*)$
of free modules over the Novikov ring
$\Lambda_{[\omega]}$, which is actually a field, see section~1.4.

We introduce a numerical invariant $\rho(\omega,g) \in [0,\infty]$,
see Definition~3 in section~1.2,
conjecturally always smaller than $\infty$, and the first purpose of this
paper is to show that if $\rho(\omega,g)<\infty$
(cf Theorem~2\therosteritem3)
\roster
\item[3]
For any $x\in\Cr_q(\omega)$ and $y\in\Cr_{q-1}(\omega)$
the sum
$$
\sum_{\hat\alpha\in\Cal P(x,y)}
I_q(x,y,\hat\alpha)e^{s[\omega](\hat\alpha)}
$$
\endroster
defines a Dirichlet series which is holomorphic in the half plane
$\{s\in\C\mid\re(s)>\rho\}$.

As a consequence, we answer positively (in the case $\rho<\infty$) a question
raised by S.~P.~Novikov. An other partial  answered to this question was 
provided by A.~V.~Pazhitnov, cf section~1.4. At this point it may be 
useful to state that we believe that $\rho$ is always smaller than 
$\infty$.

We show that the collections of above numbers can be algebraically
organized to provide a cochain complex $(\bold C^*,\boldsymbol\partial^*)$
of free modules over a much smaller ring
$\Lambda'_{[\omega],\rho}\subset\Lambda_{[\omega]}$, cf
sections~1.1 and 1.4. Actually the ring
$\Lambda'_{[\omega],\rho}$ admits, for any complex number
$s$ with $\re(s)>\rho$, a surjective ring homomorphism
$\ev_s:\Lambda'_{[\omega],\rho}\to\C$.

We define a family (in the parameter $s$, $\re(s)>\rho$)
of finite dimensional cochain complexes $(C^*,\partial^*_s)$
over the field $\C$ whose component
$C^q$ is the vector space generated by
the critical points of index $q$. For any $\re(s)>\rho$ this complex is
isomorphic to the tensor product
$\C\otimes_{\Lambda'_{[\omega],\rho}}(\bold C^*,\boldsymbol\partial^*_s)$.
With respect to the canonical base of $C^q$, the boundary map
$\partial^q_s$
can be written as a matrix whose entries
$\partial^q_s(x,y)$,
$x\in\Cr_{q+1}(\omega)$, $y\in\Cr_q(\omega)$ are functions
$s\mapsto I_{q+1,s}(x,y)$ which we show are Dirichlet series obtained from
the numbers $I_{q+1}(x,y,\hat\alpha)$.
In particular the numbers
$I_q(x,y,\hat\alpha)$ are entirely determined by the restriction of these
functions to $(a,\infty)$, for any $a>\rho$.

The second purpose of this paper is to construct (using analysis=spectral
geometry) a smooth one parameter
family of cochain complexes
$\bigl(\Omega^*_\tsm(M),d^*_t\bigr)$ which carries implicitly all
information
provided by the Novikov complex.

Precisely, given a system as above $(M,\omega,g)$, with $\omega$ a Morse
form,
Theorem~3, claims that there exists
a positive real number $T$ so that for $t\geq T$ the deRham complex
$\bigl(\Omega^*(M),d^*_t:=d+t\omega\wedge\cdot\bigr)$
decomposes canonically as a direct orthogonal sum of two complexes
$\bigl(\Omega^*_\tsm(M),d^*_t\bigr)$ and
$\bigl(\Omega^*_\tla(M),d^*_t\bigr)$. The first complex has the
component
$\Omega^q_\tsm(M)$, a finite dimensional vector space of
dimension equal to the cardinality of $\Cr_q(\omega)$.
In the case of an exact form this result is due to E.~Witten.

If $(\omega,g)$ satisfies the Morse-Smale conditions, $\rho(\omega,g)<\infty$
and one gives the orientations $o$, by Theorem~4 we show that
the integration theory provides an isomorphism between
$\bigl(\Omega^*_\tsm(M),d^*_t\bigr)$ and
$\bigl(\maps(\Cr_*(\omega),\R),\partial^*_t\bigr)$ for any $t\geq T'$, where
$T'$ is some positive real number larger than $T$ and $\rho$ discussed above.
Moreover, for $t\geq T'$ we construct a base
$E_{t,x}\in\Omega^*_\tsm(M)$, $x\in\Cr(\omega)$. With respect to this base
$d^q_t$ is a matrix whose entries are exactly the functions
$t\mapsto I_{q+1,t}(x,y)$. This results reformulates and extends results of
Helffer and Sj\"ostrand, cf \cite{HeSj85}.

Consequently the family $\bigl(\Omega^*_\tsm(M),d^*_t\bigr)$ can
be viewed as an analytic
substitute of the Novikov complex. When the base
$E_{t,x}$ is available, which is the case if $(\omega,g)$
is Morse-Smale,  $\rho(\omega, g) <\infty$ and
the orientations
$o$ are provided, this complex permits the derivation of
the numbers $I_q(x,y,\hat\alpha)$.

All these results are immediate corollaries of Theorems~1--4 stated in
section~1 and of Proposition~4 in section~5, which are of independent interest
and have many other pleasant applications.

\head 1. The results\endhead

\subhead 1.1 Novikov rings\endsubhead
Let $(M^n,\omega)$ be a pair consisting of a closed connected
smooth $n$-dimensional manifold $M$ and a closed real valued 1-form
$\omega\in\Cal Z^1(M):=\{\omega\in\Omega^1(M)\mid d\omega=0\}$.
The form $\omega$ induces the homomorphism
$[\omega]:H_1(M;\Bbb Z)\to\Bbb R$ whose image is a finitely
generated free Abelian group of rank $r$.
Denote by $\Gamma:=H_1(M;\Bbb Z)/\ker([\omega])$.
The integer $r=\rank(\Gamma)$ is called {\it degree of irrationality}
of the form $\omega$. We identify $\Gamma$ to $\Bbb Z^r$ by choosing a
base $e_1,\dotsc,e_r\in\Gamma$ with $[\omega](e_i)=\kappa_i\in\Bbb R$
positive real numbers $\Bbb Q$-linearly independent.

Let $\tilde M@>\pi>>M$ be the regular $\Gamma$ covering
associated with $H_1(M;\Bbb Z)\to \Gamma$, \ie $\tilde M$ is a
connected covering, such that for one (and hence all) $\tilde m\in\tilde M$
$$
\Img\big(\pi_1(\tilde M,\tilde m)@>{\pi_*}>>\pi_1(M,m)\big)
=\ker\big(\pi_1(M,m)@>{[\omega]}>>\R\big),
$$
where $m=\pi(\tilde m)$.
The group $\Gamma$ acts freely on $\tilde M$ with quotient space $M$.

The pull back of $\omega$ on $\tilde M$ is exact, \ie $\pi^*\omega=dh$,
with $h:\tilde M\to\Bbb R$ a smooth function.
This function is unique up to an additive constant. Given
$\tilde m\in\tilde M$ there exists a unique function $h^{\tilde m}$,
so that $\pi^*\omega=dh^{\tilde m}$ and $h^{\tilde m}(\tilde m)=0$.
In particular $\omega$ induces a function
$H:\tilde M\times\tilde M\to\Bbb R$ defined by
$$
H(\tilde x,\tilde y)=h^{\tilde m}(\tilde x)-h^{\tilde m}(\tilde y),
$$
which is independent of $\tilde m$. When there is no risk of confusion
we write $h$ for any of the function $h^{\tilde m}$. Note that
$$
h(\gamma\tilde x)=h(\tilde x)+[\omega](\gamma),
\tag 1.1
$$
for all $\gamma\in\Gamma$.

S.~P.~Novikov, see \cite{N93}, has introduced the ring
$\Lambda_{[\omega]}$, consisting of functions $f:\Gamma\to\C$ 
with the property that for any $R\in\R$ the set
$$
\big\{\gamma\in\Gamma\mid f(\gamma)\neq 0,[\omega](\gamma)\leq R\big\}
$$
is finite. The product in this ring is given by convolution, \ie
$$
(f*g)(\gamma)=\sum_{\tilde\gamma\in\Gamma}
f(\tilde\gamma)g(\tilde\gamma^{-1}\gamma).
$$
Because of the Novikov condition above this sum is actually finite and $f*g$ 
is in $\Lambda_{[\omega]}$. Since $[\omega]:\Gamma\to\R$ is injective,
or equivalently the numbers $\kappa_i$ are $\Bbb Q$-linearly independent,
$\Lambda_{[\omega]}$ is actually a field, cf \cite{HoSa95}.

Each $f\in\Lambda_{[\omega]}$ defines a Dirichlet series
$$
\hat f(s)
:=\sum_{\gamma\in\Gamma}f(\gamma)e^{-s[\omega](\gamma)}
=\sum_{n_i\in\Z}f(n_1,\dotsc,n_r)e^{-s(\kappa_1n_1+\cdots+\kappa_rn_r)},
$$
where the set of numbers $\kappa_1n_1+\cdots+\kappa_rn_r$ with
$f(n_1,\dotsc,n_r)\ne 0$ is a strictly increasing sequence of real numbers
$\lambda_1<\lambda_2<\cdots$ which is either finite or is tending to 
$+\infty$.

Recall that if $(\lambda_n)_{n\in\N}$ is an increasing sequence of real
numbers tending to $+\infty$ a Dirichlet series
with exponents $\lambda_n$ is a series of the form
$\hat f:=\sum a_ne^{-s\lambda_n}$,
$a_n\in\C$, $s\in\C$. If the series converges for $s_0$, it defines
a holomorphic function on the open half plane
$\{s\in\C\mid\re(s)>\re(s_0)\}$ so there exists
$\rho(\hat f)\in\R\cup\{\infty\}$,
referred to as the abscissa of convergence of
$\hat f$, making the series a holomorphic function on
$\{s\in\C\mid\re(s)>\rho\}$.
Note that $\widehat{f*g}=\hat f\cdot\hat g$ and
$\rho(\widehat{f*g})\leq\sup\{\rho(\hat f),\rho(\hat g)\}$.

Let $\Lambda'_{[\omega],\rho}$ be the the subring of $\Lambda_{[\omega]}$
consisting of elements $f\in\Lambda_{[\omega]}$ whose corresponding
Dirichlet series is convergent for any $s\in\C$ with $\re(s)>\rho$. 
Any such $s$ gives rise to an evaluation homomorphism
$$
\ev_s:\Lambda'_{[\omega],\rho}\to\C,\quad
f\mapsto\hat f(s),
\tag 1.2
$$
and let $\iota:\Lambda'_{[\omega],\rho}\to\Cal F$ be the obvious ring 
homomorphism obtained by restricting the holomorphic function defined by the
element in $\Lambda'_{[\omega],\rho}$ to the interval $(\rho,\infty)$, and
where $\Cal F$ denotes the ring of germs at $+\infty$ of $\C$-valued
smooth functions $f:(a,\infty)\to\C$, where $a\in\R$. Clearly $\iota$ is
injective and the general theory of Dirichlet series (or almost periodic
functions) permits to recover the coefficients $f(n_1,\dotsc,n_r)$ from
the germ $\iota(f)$, cf \cite{Se73}.

\subhead 1.2 Morse-Smale condition and the invariant $\rho$\endsubhead
Recall that for $x\in\Cr(\omega):=\Zeros(\omega)$ the Hessian of
$\omega$ at $x$ is
$$
H_x\omega:T_xM\times T_xM\to\R,\quad
(H_x\omega)(X,Y):=(\nabla_X\omega)(Y),
\tag 1.3
$$
where $\nabla$ is any linear connection on $M$. $H_x\omega$ does not depend
on the connection and is symmetric since $\omega$ is closed. The closed
1-form $\omega$ is called Morse form if $H_x\omega$ is non-degenerate
for every $x\in\Cr(\omega)$. The index of $\omega$ at $x\in\Cr(\omega)$
is the index of $H_x\omega$. By the Morse lemma, for any $x\in\Cr(\omega)$ 
there exists an open neighborhood
$U_x$ of $x$, positive real numbers $c_x,\epsilon_x$ and a diffeomorphism
$\theta_x:(U_x, x)\to(D^n(\epsilon_x),0)$, where  $D^n(r)$ denotes the
open disc of radius $r$ in $\R^n$ centered at $0$, so that
$$
(\theta_x^{-1})^*\omega
=d\big(-c_x(x^2_1+\cdots+x^2_k)+c_x(x^2_{k+1}+\cdots+x^2_n)\big),
\tag a
$$
where $k=\index(x)$.
In what follows we consider systems $(M,\omega,g)$ where $M$ is a closed
manifold, $\omega$ a closed 1-form as above and $g$ is a Riemannian metric.

Let $\grad_g\omega$ be the unique vector field which corresponds to
$\omega$ by the bijective correspondence between vector fields and closed
1-forms provided by the Riemannian metric $g$ and set $X:=-\grad_g\omega$.
For each $x\in M$ denote by $\gamma_x(t)$ the trajectory of $X$ with
$\gamma_x(0)=x$. For $x\in\Cr(\omega)$ denote by $W^\pm_x$ the sets
$$
W^\pm_x=\bigl\{y\bigm|\lim_{t\to\pm\infty}\gamma_y(t)=x\bigr\}.
$$
They will be referred to as the stable \resp unstable sets of the
critical point $x$.

\definition{Definition 1 (Morse pairs)}
The pair $(\omega,g)$ is called a Morse pair if for any
$x\in\Cr(\omega)$ there exists $\epsilon_x,c_x$ and $\theta_x$, so that
\thetag{a} and the following condition \thetag{b} are satisfied.
$$
(\theta_x^{-1})^*g=dx^1\otimes dx^1+\cdots+dx^n\otimes dx^n
\tag b
$$
\enddefinition

In view  of the theorem of existence, uniqueness and smooth dependence on the
initial conditions for the solutions of ordinary differential equations,
the fact that $(\omega,g)$ is a Morse pair implies that
$W_x^-$ \resp $W_x^+$ is
the image by a smooth one to one immersion $i^-_x:\R^k\to M$ \resp
$i^+_x:\R^{n-k}\to M$, where $k=\index(x)$.
Denote by $h^x:\R^k\to\R$ the unique smooth map which satisfies
$(i^-_x)^*\omega=dh^x$, $h^x(0)=0$ and by $g^x:=(i^-_x)^*g$ the pull back
of the Riemannian metric $g$ by the immersion $i^-_x$, which is a
Riemannian metric on $\R^k$. 

\definition{Definition 2 (Morse-Smale condition)}
The pair $(\omega,g)$ is called Morse-Smale if it is a Morse pair and in
addition for any $x,y\in\Cr(\omega)$, $i_x^- $ and $i_y^+$ are transversal.
Note that if $(\omega,g)$ is a Morse \resp Morse-Smale pair then so is
$(t\omega,g)$ for any $0\neq t\in\R$.
\enddefinition

Denote by  $\Cal G$ be the set of smooth Riemannian metrics on $M$.
Let $U\subset M$ be open and $g\in\Cal G$. Denote by $\Cal G_{g,U}$
the set
$$
\Cal G_{g,U}:=\bigl\{g'\in\Cal G\bigm|
\forall x\in M\setminus U:g'(x)=g(x)\bigr\}.
$$
The following almost obvious result establishes the existence of Morse
forms and Morse pairs.

\proclaim{Proposition 1}
Suppose $M$ is a closed manifold. Then the following holds:
\roster
\item
The set of Morse forms is open and dense subset of $\Cal Z^1(M)$
equipped with the $C^1$-topology.
\item
Let $\omega$ be a Morse form, $g$ a Riemannian metric and let $U$ be a
neighborhood of $\Cr(\omega)$. Then the set of metrics
$g'\in\Cal G_{g,U}$, so that $(\omega,g')$ is a Morse pair is dense in
$\Cal G_{g,U}$ with respect to the $C^0$-topology.
\endroster
\endproclaim

The following proposition establishes the existence of Morse-Smale pairs.
Its proof can be derived from Kupka-Smale's theorem, cf \cite{Pe67}.
In section~2 we will give an alternative proof on the lines of \cite{Sch93}.

\proclaim{Proposition 2}
Let $(\omega,g)$ be a Morse pair, $\varepsilon>0$ small and set
$$
U:=\bigcup_{z\in\Cr(\omega)}B(z,\varepsilon)\setminus
\overline{B(z,\tfrac\varepsilon2)}.
$$
Then there exists a Banach manifold $G\subseteq\Cal G_{g,U}$
of smooth Riemannian metrics, which is dense in $\Cal G_{g,U}$ with respect to
the $L^2$-topology, and a residual subset
$G'\subset G$, such that for any $g'\in G'$ the pair
$(\omega,g')$ is Morse-Smale.
\endproclaim

\definition{Definition 3 (The invariant $\rho$)}
For a Morse-Smale pair $(\omega,g)$ we denote
$$
\rho(\omega,g)
:=\inf\bigl\{a\in\R_+\bigm|
\forall x\in\Cr(\omega):\int_{\R^{\index(x)}}e^{ah^x}\Vol_{g^x}<\infty\bigr\}.
$$
It is conceivable that there are no such positive real numbers $a$,
in which case we put  $\rho(\omega,g)=\infty$.
\enddefinition

We believe that always $\rho(\omega,g)<\infty$. There are plenty of
examples where $\rho(\omega,g)=0$.
\footnote{
If $\dim(M)\leq 2$ and $(\omega,g)$ is a Morse-Smale pair then
$\rho(\omega,g)=0$. Indeed, for $\index(x)=1$ this follows from Lemma~3 in
section~3, below. Moreover for $\index(x)=n=\dim(M)$ one has for all $a\geq 0$
$$
\int_{\R^n}e^{ah^x}\Vol_{g^x}\leq\Vol(M)<\infty,
$$
since $i_x^-:\R^n\to M$ is a one to one immersion.
}
\footnote{
If $(M_i,\omega_i,g_i)$, $i=1,2$ are Morse-Smale pairs
then $\bigl(M_1\times M_2,\pi_1^*\omega_1+\pi_2^*\omega_2,
\pi_1^*g_1+\pi_2^*g_2\bigr)$ is a Morse-Smale
pair and
$\rho\bigl(\pi_1^*\omega_1+\pi_2^*\omega_2,\pi_1^*g_1+\pi_2^*g_2\bigr)\leq
\sup\bigl\{\rho(\omega_1,g_1),\rho(\omega_2,g_2)\bigr\}$, where
$\pi_i:M_1\times M_2\to M_i$ denotes the canonical projection.
}
The invariant $\rho(\omega,g)$ will be discussed in a forthcoming paper.

\subhead 1.3 Compactification\endsubhead
Let $(\omega,g)$ be a Morse-Smale pair. Denote the set of critical points
of $h$ by $\Cr(h):=\Cr(dh)=\pi^{-1}(\Cr(\omega))$. Recall that for
$\tilde x\in\Cr(h)$ one has the stable and unstable manifolds
$W^\pm_{\tilde x}$ of the negative gradient flow of $h$. One can also
consider the immersions $i^+_{\tilde x}:\R^{n-k}\to\tilde M$ and
$i^-_{\tilde x}:\R^k\to\tilde M$, where $k=\index(\tilde x)$, which in
this case are embeddings. The submanifolds $W^\pm_{\tilde x}$ are exactly
their images. Note that with the notation $h^x$ introduced in section~1.2, 
one has $h^x=h^{\tilde x}\circ i^-_{\tilde x}$, for any $\tilde x$ with
$\pi(\tilde x)=x$.

\remark{Observation 1}
$\Cr(h)\subseteq\tilde M$ is a discrete subset. $\Gamma$ acts transitively
and freely on $\Cr(h)$ with quotient set $\Cr(\omega)$.
\endremark

\remark{Observation 2}
$\pi:h^{-1}(c)\to M$ is injective, for all $c\in\R$. In particular any
critical level of $h$ contains only finitely many critical points.
\endremark

\remark{Observation 3}
If $\tilde x\in\Cr(h)$, then
$\pi:W^\pm_{\tilde x}\to W^\pm_{\pi(x)}\subset M$ is an injective immersion.
The Morse-Smale condition from Definition~2 is equivalent to the
transversality of $W^-_{\tilde x}$ and $W^+_{\tilde y}$ for any two
$\tilde x,\tilde y\in\Cr(h)$.
\endremark

\remark{Observation 4}
The Morse-Smale condition implies that $\Cal M(\tilde x,\tilde y):=
W^-_{\tilde x}\cap W^+_{\tilde y}$ is a submanifold of $\tilde M$ of 
dimension $\index(x)-\index(y)$. The manifold
$\Cal M(\tilde x,\tilde y)$ is equipped with the action
$\mu:\R\times\Cal M(\tilde x,\tilde y)\to\Cal M(\tilde x,\tilde y)$,
defined by $\mu(t, z)=\gamma_z(t)$. If $\tilde x\neq\tilde y$
the action $\mu$ is free and we denote the quotient
$\Cal M(\tilde x,\tilde y)/\R$ by $\Cal T(\tilde x,\tilde y)$.
$\Cal T(\tilde x,\tilde y)$ is a smooth manifold of dimension
$\index(x)-\index(y)-1$, possibly empty, diffeomorphic to the submanifold
$h^{-1}(c)\cap\Cal M(\tilde x,\tilde y)$, where $c$ is any regular
value of $h$ with $h(\tilde x)>c>h(\tilde y)$. Note that if
$\index(\tilde x)\leq\index(\tilde y)$, and $\tilde x\neq\tilde y$,
in view of the transversality requested by the Morse-Smale condition,
$\Cal M(\tilde x,\tilde y)=\emptyset$. If $\tilde x=\tilde y$, then
$W_{\tilde x}^-\cap W_{\tilde x}^+=\tilde x$. The elements of
$\Cal T(\tilde x,\tilde y)$ will be referred to as
the unparameterized trajectories from $\tilde x$ to $\tilde y$.
\endremark

\definition{Definition 4 (Broken trajectories)}
An unparameterized broken trajectory from $\tilde x\in\Cr(h)$ to
$\tilde y\in\Cr(h)$ is an element of
$$
\Cal B(\tilde x,\tilde y):=\bigcup_{
\Sb k\geq 0,\ \tilde y_0,\dotsc,\tilde y_{k+1}\in\Cr(h)\\
\tilde y_0=\tilde x,\ \tilde y_{k+1}=\tilde y\\
\index(\tilde y_i)>\index(\tilde y_{i+1})
\endSb}
\Cal T(\tilde y_0,\tilde y_1)
\times\cdots\times\Cal T(\tilde y_k,\tilde y_{k+1}).
$$
An unparameterized broken trajectory from $\tilde x\in\Cr(h)$ to the level
$\lambda\in\R$ is an element of
$$
\Cal B(\tilde x;\lambda):=
\bigcup_{
\Sb k\geq 0,\ \tilde y_0,\dotsc,\tilde y_k\in\Cr(h)\\
\tilde y_0=\tilde x\\
\index(\tilde y_i)>\index(\tilde y_{i+1})
\endSb}
\Cal T(\tilde y_0,\tilde y_1)\times\cdots\times
\Cal T(\tilde y_{k-1},\tilde y_k)\times
(W^-_{\tilde y_k}\cap h^{-1}(\lambda)).
$$
\enddefinition

Clearly, if $\lambda>h(\tilde x)$ then
$\Cal B(\tilde x;\lambda)=\emptyset$. There is an obvious way to
regard $\Cal B(\tilde x,\tilde y)$ \resp $\Cal B(\tilde x;\lambda)$
as a subset of $C^0\big([h(\tilde y),h(\tilde x)],\tilde M\big)$ \resp
$C^0\big([\lambda,h(\tilde x)],\tilde M\big)$, by parameterizing a broken
trajectory by the value of $h$. This leads to the following
characterization and implicitly to a canonical
parameterization of an unparameterized broken trajectory.

\remark{Observation 5}
Let $\tilde x,\tilde y\in\Cr(h)$ and set $a:=h(\tilde y)$, $b:=h(\tilde x)$.
The parameterization above defines a one to one correspondence between
$\Cal B(\tilde x,\tilde y)$ and the set of continuous mappings
$\gamma:[a,b]\to\tilde M$, which satisfy the following two properties:
\roster
\item
$h(\gamma(s))=a+b-s$, $\gamma(a)=\tilde x$ and $\gamma(b)=\tilde{y}$.
\item
There exists a finite collection of real numbers
$a=s_0<s_1<\cdots<s_{r-1}<s_r=b$, so that $\gamma(s_i)\in\Cr(h)$
and $\gamma$ restricted to $(s_i,s_{i+1})$ has derivative at any
point in the interval $(s_i,s_{i+1})$, and the derivative satisfies
$$
\gamma'(s)=\frac{-\grad_gh}{\|\grad_gh\|^2}\big(\gamma(s)\big).
\tag 1.4
$$
\endroster
Similarly the elements of $\Cal B(\tilde x;\lambda)$ correspond to
continuous mappings $\gamma:[\lambda,b]\to\tilde M$, which satisfies
\therosteritem1 and \therosteritem2, with $a$ replaced by $\lambda$.
\endremark

In section~3 we will verify the following

\proclaim{Proposition 3}
Let $(\omega,g)$ be a Morse-Smale pair, $\tilde x,\tilde y\in\Cr(h)$ and
$\lambda\in\R$. Then:
\roster
\item
$\Cal B(\tilde x,\tilde y)$ is compact, with the topology induced from
$C^0\big([h(\tilde y),h(\tilde x)],\tilde M\big)$.
\item
$\Cal B(\tilde x;\lambda)$ is compact, with the topology induced from
$C^0\big([\lambda,h(\tilde x)],\tilde M\big)$.
\endroster
\endproclaim

For $\tilde y_0,\dotsc,\tilde y_k\in\Cr(h)$ with
$\index(\tilde y_i)>\index(\tilde y_{i+1})$, consider the smooth map
$$
i_{\tilde y_0,\dotsc,\tilde y_k}:
\Cal T(\tilde y_0,\tilde y_1)\times\cdots\times
\Cal T(\tilde y_{k-1},\tilde y_k)\times W_{\tilde y_k}^-
\to\tilde M,
$$
defined by $i_{\tilde y_0,\dotsc,\tilde y_k}
(\gamma_1,\dotsc,\gamma_k,\tilde y):=i_{\tilde y_k}(\tilde y)$, where
$i_{\tilde x}:W_{\tilde x}^-\to\tilde M$ denotes the inclusion.

\definition{Definition 5 (Completed unstable manifold)}
For $\tilde x\in\Cr(h)$ define
$$
\hat W^-_{\tilde x}:=\bigcup_{
\Sb k\geq 0,\ \tilde y_0,\dotsc,\tilde y_k\in\Cr(h)\\
\tilde y_0=\tilde x\\
\index(\tilde y_i)>\index(\tilde y_{i+1})
\endSb}
\Cal T(\tilde y_0,\tilde y_1)\times\cdots\times
\Cal T(\tilde y_{k-1},\tilde y_k)\times W^-_{\tilde y_k}.
$$
Moreover, let $\hat i_{\tilde x}:\hat W^-_{\tilde x}\to\tilde M$ denote
the mapping, whose restriction to $\Cal T(\tilde y_0,\tilde y_1)
\times\cdots\times\Cal T(\tilde y_{k-1},\tilde y_k)\times W^-_{\tilde y_k}$
is given by $i_{\tilde y_0,\dotsc,\tilde y_k}$ and let
$\hat h_{\tilde x}:=h\circ\hat i_{\tilde x}:\hat W^-_{\tilde x}\to\R$.
\enddefinition

To formulate the next result we need an  additional concept, smooth
manifold with corners. Recall that an $n$-dimensional manifold $P$
with corners is a paracompact Hausdorff space equipped with a
maximal smooth atlas with charts $\varphi:U\to\varphi(U)\subseteq\R^n_+$,
where $\R^n_+=\{(x_1,x_2,\dotsc,x_n)\mid x_i\geq 0\}$. The collection of
points of $P$ which correspond (by some and then by any chart) to points
in $\R^n$ with exactly $k$ coordinates equal to zero is a well defined
subset of $P$ and it will be denoted by $P_k$. It has a structure of a
smooth $(n-k)$-dimensional manifold.
$\partial P=P_1\cup P_2\cup\cdots\cup P_n$ is a closed subset which
is a topological manifold, and $(P,\partial P)$ is a topological
manifold with boundary $\partial P$.

\proclaim{Theorem 1}
Let $(\omega,g)$ be a Morse-Smale pair.
\roster
\item
For any two critical points $\tilde x,\tilde y\in\Cr(h)$ the
smooth manifold $\Cal T(\tilde x,\tilde y)$ has
${\Cal B}(\tilde x,\tilde y)$ as a canonical
compactification. Moreover ${\Cal B}(\tilde x,\tilde y)$ has the structure
of a compact smooth manifold with corners, and
$$
{\Cal B}(\tilde x,\tilde y)_k
=\bigcup_{
\Sb\tilde y_0,\dotsc,\tilde y_{k+1}\in\Cr(h)\\
\tilde y_0=\tilde x,\ \tilde y_{k+1}=\tilde y\\
\index(\tilde y_i)>\index(\tilde y_{i+1})\endSb}
\Cal T(\tilde y_0,\tilde y_1)
\times\cdots\times\Cal T(\tilde y_k,\tilde y_{k+1}),
$$
especially $\Cal B(\tilde x,\tilde y)_0=\Cal T(\tilde x,\tilde y)$.
\item
For any critical point $\tilde x\in\Cr(h)$, $\hat W_{\tilde x}^-$ has a
canonical structure of a smooth manifold with corners, and
$$
(\hat W_{\tilde x}^-)_k
=\bigcup_{
\Sb\tilde y_0,\dotsc,\tilde y_k\in\Cr(h)\\
\tilde y_0=\tilde x\\
\index(\tilde y_i)>\index(\tilde y_{i+1})
\endSb}
\Cal T(\tilde y_0,\tilde y_1)\times\cdots\times
\Cal T(\tilde y_{k-1},\tilde y_k)\times W_{\tilde y_k}^-,
$$
especially $(\hat W_{\tilde x}^-)_0=W_{\tilde x}^-$. Moreover
$\hat i_{\tilde x}$ and $\hat h_{\tilde x}$ are smooth and proper maps,
and $\hat i_{\tilde x}$ is a closed map.
\endroster
\endproclaim

The proof of this theorem will be given in section~4.
Propositions~1--3 are known in literature and Theorem~1,
can be also found in \cite{L95}. Our proof of Theorem~1 is however
different from the one sketched in \cite{L95} and we hope more conceptual.
It also has the virtue that it extends essentially word by word to
Bott-Smale pairs.

\subhead 1.4 Novikov complexes\endsubhead
Let $(\omega,g)$ be a Morse-Smale pair. For
any $x\in\Cr(\omega)$ choose an orientation $\Cal O_x$ in $W^-_x$ and
denote the collection of these orientations by $o$. Via $\pi$ these
orientations induce orientations $\Cal O_{\tilde x}$ on $W^-_{\tilde x}$.
Denote by $\Cal X_q:=\Cr_q(h)=\{\tilde x\in\Cr(h)\mid\index(x)=q\}$.
Theorem~1 implies the existence of the map
$$
I_q:\Cal X_q\times\Cal X_{q-1}\to\Z
$$
defined as follows:

If $\Cal T(\tilde x,\tilde y)=\emptyset$ put
$I_q(\tilde x,\tilde y)=0$. If $\Cal T (\tilde x,\tilde y)\neq\emptyset$,
then for any $\gamma\in\Cal T(\tilde x,\tilde y)$ the set
$\gamma\times W_{\tilde y}^-$ appears as an open set of the boundary
$\partial\hat W_{\tilde x}^-$ and the orientation $\Cal O_{\tilde x}$
induces an orientation on it. If this is the same as the
orientation $\Cal O_{\tilde y}$, we set  $\varepsilon(\gamma)=+1$,
otherwise we set $\varepsilon(\gamma)=-1$. Now define
$I_q(\tilde x,\tilde y)$ by
$$
I_q(\tilde x,\tilde y)
:=\sum_{\gamma\in\Cal T(\tilde x,\tilde y)}\varepsilon(\gamma),
$$
which is a finite sum by Proposition~3\therosteritem1. The following result
establishes the main properties of the numbers $I_q(\tilde x,\tilde y)$.

\proclaim{Theorem 2}
Suppose $(\omega,g)$ is a Morse-Smale pair. Then:
\roster
\item
$I_q(\gamma\tilde x,\gamma\tilde y)=I_q(\tilde x,\tilde y)$, for all
$\gamma\in\Gamma$.
\item
For all $\tilde x\in\Cal X_q$ and $\tilde y\in\Cal X_{q-2}$
the sum below contains only finitely many nonzero terms and one has
$$
\sum_{\tilde z\in\Cal X_{q-1}}
I_q(\tilde x,\tilde z)I_{q-1}(\tilde z,\tilde y)=0.
\tag 1.5
$$
\item
For any $\tilde x\in\Cal X_q$, $\tilde y\in\Cal X_{q-1}$
and any $s\in\C$ with $\re(s)>\rho(\omega,g)$ the sum
$$
\sum_{\gamma\in\Gamma}I_q(\gamma\tilde x,\tilde y)e^{-s[\omega](\gamma)}
\tag 1.6
$$
is convergent.
\endroster
\endproclaim

If $\rho(\omega,g)<\infty$ Theorem~2\therosteritem3 provides a positive
answer to the following conjecture formulated by S.~P.~Novikov, cf \cite{N93},
\cite{A90} and \cite{Pa98}.

\proclaim{Conjecture (Novikov)}
For any Morse-Smale pair $(\omega,g)$ and any two critical points
$\tilde x\in\Cal X_q$ and $\tilde y\in\Cal X_{q-1}$ the integers
$I_q(\gamma\tilde x,\tilde y)$, $\gamma\in\Gamma\cong\Z^r$
have at most exponential growth. More precisely, there exist constants
$C_{\tilde x,\tilde y}, M_{\tilde x,\tilde y}\in\R$, such that
$$
I_q(\tilde\gamma\tilde x,\tilde y)\leq
C_{\tilde x,\tilde y}e^{M_{\tilde x,\tilde y}[\omega](\gamma)},
$$
for all $\gamma\in\Gamma$.
\endproclaim

In \cite{Pa98} A.~V.~Pazhitnov has verified this conjecture for a generic
subset of the set of Riemannian metrics $g$, for which $(\omega,g)$ is
Morse-Smale.

\remark{Remark}
Note, that the numbers $I_q(\tilde x,\tilde y)$ can be defined without any
reference to the covering $\pi:\tilde M\to M$.
Indeed, let $x\in M$ and choose $\tilde x\in\tilde M$, such that
$\pi(\tilde x)=x$. Then there exists a natural one to one correspondence
between $\Cal P(x,y)$ and $\pi^{-1}(y)$, given by lifting a path $\alpha$
from $x$ to $y$, to a path starting at $\tilde x$ and looking at the
endpoint, which does only depend on $\hat\alpha$. In view of
Theorem~2\therosteritem1,
$I_q(x,y,\hat\alpha):=I_q(\tilde x,\tilde y)$ is well defined, \ie
independent of the choice of $\tilde x$, and the
formulas \thetag{1.5} and \thetag{1.6} become \therosteritem2 and
\therosteritem3 in the introduction.
\endremark

\definition{Definition 6 (Novikov condition)}
We say a map $f:\Cal X_q\to\C$ has property
\roster\widestnumber\item{(N$_\rho$)}
\item"(N)"
if for any $\tilde x\in\Cal X_q$ and any $R\in\R$ the set
$\{\gamma\in\Gamma\mid f(\gamma\tilde x)\neq 0,[\omega](\gamma)\leq R\}$
is finite, and we say it has property
\item"(N$_\rho$)"
if for any $\tilde x\in\Cal X_q$ and $s\in\C$ with $\re(s)>\rho$ the series
$\sum_\gamma f(\gamma\tilde x)e^{-s[\omega](\gamma)}$ is convergent.
\endroster
\enddefinition

Let $N\bold C^q$ denote the $\C$-vector space of functions $\Cal X_q\to\C$,
which satisfy (N). For $\lambda\in\Lambda_{[\omega]}$ and
$f\in N\bold C^q$ we set
$$
(\lambda*f)(\tilde x)
:=\sum_{\gamma\in\Gamma}\lambda(\gamma)f(\gamma^{-1}\tilde x).
$$
In this way $N\bold C^q$ becomes a free $\Lambda_{[\omega]}$-module of
finite rank equal to the cardinality $\Cr_q(\omega)$. Moreover let
$\bold C^q$ denote the subspace of functions, which satisfy
(N) and (N$_\rho$). The formula above also makes $\bold C^q$ a free
$\Lambda'_{[\omega],\rho}$-module of the same rank as $N\bold C^q$. Note, that
every section $\sigma:\Cr_q(\omega)\to\Cal X_q$, \ie $\pi\circ\sigma=\id$,
defines a base for both $N\bold C^q$ and $\bold C^q$, namely
$\{\delta_{\sigma(x)}\mid x\in\Cr_q(\omega)\}$, where
$\delta_{\tilde x}:\Cr(h)\to\C$ is the Kronecker function,
$\delta_{\tilde x}(\tilde y)=\delta_{\tilde x,\tilde y}$.

For $\tilde y\in\Cal X_q$ we define $\boldsymbol\partial^q(\delta_{\tilde y})\in
\maps(\Cal X_{q+1},\C)$ by $(\boldsymbol\partial^q(\delta_{\tilde y}))(\tilde x)
:=I_{q+1}(\tilde x,\tilde y)$. Theorem~2\therosteritem3 shows that
$\boldsymbol\partial^q(\delta_{\tilde y})$ satisfies (N$_\rho$) and
Corollary~1 in section~3 shows, that it also satisfies (N),
\ie $\boldsymbol\partial^q(\delta_{\tilde y})\in\bold C^{q+1}\subseteq
N\bold C^{q+1}$. From Theorem~2\therosteritem1 one
gets $\boldsymbol\partial^q(\delta_{\gamma\tilde y})
=\delta_\gamma*\boldsymbol\partial^q(\delta_{\tilde y})$. This equivariance
property and the fact that $N\bold C^q$ and $\bold C^q$ are free modules
shows, that $\boldsymbol\partial^q$ extends uniquely to a
$\Lambda_{[\omega]}$ \resp $\Lambda'_{[\omega],\rho}$-linear map
$$
\boldsymbol\partial^q:N\bold C^q\to N\bold C^{q+1}\quad\text{\resp}\quad
\boldsymbol\partial^q:\bold C^q\to\bold C^{q+1},
$$
both given by the formula
$$
\boldsymbol\partial^q(f)(\tilde x)
=\sum_{\tilde y\in\Cal X_q}I_{q+1}(\tilde x,\tilde y)f(\tilde y).
$$
Theorem~2\therosteritem2 immediately shows
$\boldsymbol\partial^{q+1}\circ\boldsymbol\partial^q=0$, by
checking
it on the elements $\delta_{\tilde y}$. So we have two cochain complexes,
the Novikov complex $(N\bold C^*,\boldsymbol\partial^*)$ and
$(\bold C^*,\boldsymbol\partial^*)$ and a natural isomorphism
$$
\Lambda_{[\omega]}\otimes_{\Lambda'_{[\omega],\rho}}(\bold C^*,\boldsymbol\partial^*)
\cong
(N\bold C^*,\boldsymbol\partial^*),\quad\lambda\otimes f\mapsto\lambda*f.
$$

Define $C^q:=\maps(\Cr_q(\omega),\C)$ and note that this is a finite
dimensional $\C$-vector space.
Suppose $\rho:=\rho(\omega,g)<\infty$ and let $s\in\C$ with $\re(s)>\rho$. For
$x\in\Cr_q(\omega)$ and $y\in\Cr_{q-1}(\omega)$ choose
$\tilde y\in\Cal X_{q-1}$, such that $\pi(\tilde y)=y$ and define
$$
I_{q,s}(x,y):=\sum_{\tilde x\in\pi^{-1}(x)}
I_q(\tilde x,\tilde y)e^{-sH(\tilde x,\tilde y)}\in\C,
$$
which converges by Theorem~2\therosteritem3 and does not depend on the
choice of $\tilde y$. The map $H(\tilde x,\tilde y)
:=h(\tilde x)-h(\tilde y)$ was introduced in section~1.1.
Moreover it follows from Theorem~2\therosteritem2, that one has
$$
\sum_{z\in\Cr_{q-1}(\omega)}I_{q,s}(x,z)I_{q-1,s}(z,y)=0,
$$
for all $x\in\Cr_q(\omega)$ and $y\in\Cr_{q-2}(\omega)$. So for every
$s\in\C$ with $\re(s)>\rho$ we get another cochain complex $(C^*,\partial^*_s)$,
where
$$
\partial^q_s(f)(x):=\sum_{y\in\Cr_q(\omega)}I_{q+1,s}(x,y)f(y).
$$

Next define an evaluation map $\ev^h_s:\bold C^*\to C^*$, by
$$
\ev^h_s(f)(x):=\sum_{\tilde x\in\pi^{-1}(x)}
f(\tilde x)e^{-sh(\tilde x)}.
$$
This depends on the choice of $h$, but if one changes $h$ it changes
only by a nonzero multiplicative constant in $\C$. One easily checks
$\ev^h_s\circ\boldsymbol\partial^q=\partial^q_s\circ\ev^h_s$, \ie
$$
\ev^h_s:(\bold C^*,\boldsymbol\partial^*)\to(C^*,\partial^*_s)
$$
is a chain mapping. Moreover one has
$\ev^h_s(\lambda*f)=\ev_s(\lambda)\cdot\ev^h_s(f)$, where
$\ev_s:\Lambda'_{[\omega],\rho}\to\C$ is the evaluation map from \thetag{1.2}.
Therefore
$$
\C\otimes_{\Lambda'_{[\omega],\rho}}(\bold C^*,\boldsymbol\partial^*)\cong
(C^*,\partial^*_s),\quad
z\otimes f\mapsto z\cdot\ev^h_s(f)
$$
is an isomorphism of cochain complexes over $\C$. Here the
$\Lambda'_{[\omega],\rho}$-module structure on $\C$ is the one given by
$\ev_s:\Lambda'_{[\omega],\rho}\to\C$.

Finally let $\Omega^*(M;\C):=\Omega^*(M)\otimes\C$ denote the
$\C$-valued differential forms on $M$ and consider
$$
d^q_s:\Omega^q(M;\C)\to\Omega^{q+1}(M;\C),\quad
d^q_s(\alpha):=d\alpha+s\omega\wedge\alpha.
$$
Since $\omega$ is closed
one has $d^{q+1}_s\circ d^q_s=0$. So
$\big(\Omega^*(M;\C),d^*_s\big)$
is a cochain complex and for $\re(s)>\rho$ one has a chain mapping
$$
\Int_s:\big(\Omega^*(M;\C),d^*_s\big)\to(C^*,\partial^*_s),\quad
\Int_s(\alpha)(x):=\int_{W^-_{\tilde x}}e^{sh^{\tilde x}}\pi^*\alpha,
$$
where $\tilde x\in\Cr(h)$, such that $\pi(\tilde x)=x$, and $h^{\tilde x}$
is the unique $h$, such that $h^{\tilde x}(\tilde x)=0$. The integral
converges because of Proposition~4\therosteritem1 and is obviously
independent of the choice of $\tilde x$. Proposition~4\therosteritem2 shows
that $\Int_s$ intertwines
the differentials. Proposition~4 is stated in section~5 below. Theorem~4 in
the next section implies, that $\Int_s$
induces an isomorphism in cohomology.

The cochain complex $(C^*,\partial^*_s)$ can be regarded as a smooth
family of cochain complexes of finite dimensional vector spaces which,
in view of the fact that the cohomology of $(C^*,\partial^*_s)$ 
does not change dimension for large $s$, is a smooth bundle of cochain 
complexes, for large $s$.

\subhead 1.5 Witten-Helffer-Sj\"ostrand theory\endsubhead
Let $M$ be a closed manifold and $\omega$
a closed 1-form. For $t\in\R$ consider the complex
$\bigl(\Omega^*(M),d^*_t\bigr)$ with differential
$$
d^q_t:\Omega^q(M)\to\Omega^{q+1}(M),\quad
d^q_t(\alpha):=d\alpha+t\omega\wedge\alpha.
$$
Clearly $d^q_0=d^q$.

Recall, that on an oriented $n$-dimensional Riemannian manifold $(M,g)$
one has the Hodge-star operator $*:\Omega^q(M)\to\Omega^{n-q}(M)$. It
is a zero order operator and satisfies
$$
*\circ*=(-1)^{q(n-q)}\id:\Omega^q(M)\to\Omega^q(M).
$$
One defines the fiberwise scalar product
$$
\llangle\cdot,\cdot\rrangle:\Omega^q(M)\times\Omega^q(M)\to\Omega^0(M),
\quad\llangle\alpha_1,\alpha_2\rrangle:=*^{-1}(\alpha_1\wedge*\alpha_2)
$$
and the formal adjoint of $d^q_t$,
$(d^q_t)^\sharp:\Omega^{q+1}(M)\to\Omega^q(M)$,
$$
(d^q_t)^\sharp(\alpha)
=(-1)^{nq+1}*d^{n-q-1}_t(*\alpha)
=(d^q)^\sharp(\alpha)+i_{\grad_g\omega}\alpha.
$$

The fiberwise scalar products $\llangle\cdot,\cdot\rrangle$ and the
operators $(d^q_t)^\sharp$ are independent of the orientation of $M$.
They can even be defined (first locally and
then, being differential operators, globally) for an arbitrary
Riemannian manifold, not necessarily orientable.
Moreover one has the scalar product
$$
\Omega^q(M)\times\Omega^q(M)\to\R,\quad
\langle\alpha_1,\alpha_2\rangle
:=\int_M\alpha_1\wedge*\alpha_2
=\int_M\llangle\alpha_1,\alpha_2\rrangle\Vol.
$$
The operators $(d^q_t)^\sharp$ are formal adjoints of $d^q_t$, more
precisely
$$
\bigl\langle d^q_t(\alpha_1),\alpha_2\bigr\rangle
=\bigl\langle\alpha_1,(d^q_t)^\sharp(\alpha_2)\bigr\rangle.
$$
Next we introduce the Witten Laplacian for the closed 1-form $\omega$,
$$
\Delta^q_t:\Omega^q(M)\to\Omega^q(M),\quad
\Delta^q_t(\alpha):=(d^q_t)^\sharp\bigl(d^q_t(\alpha)\bigr)
+d^{q-1}_t\bigl((d^{q-1}_t)^\sharp(\alpha)\bigr).
$$
This is a second order differential operator, and
$\Delta^q_0=\Delta^q$, the Laplace-Beltramy operator. The operators
$\Delta^q_t$ are elliptic, selfadjoint
and positive, hence their spectra
$\spect(\Delta^q_t)$ lie on
$[0,\infty)$. Finally one has
$$
\ker(\Delta^q_t)=\bigl\{\alpha\in\Omega^q(M)\bigm| d^q_t(\alpha)=0,
(d^{q-1}_t)^\sharp (\alpha)=0\bigr\}.
$$

The following result extends a result due to E.~Witten in the case that
$\omega$ is exact.

\proclaim{Theorem 3}
Suppose that $(\omega,g)$ is a Morse pair.
There exist the
constants
$C_1$, $C_2$, $C_3$ and $T_0$ depending on $(\omega,g)$, so that for
any
$t\geq T_0$
\roster
\item
$\spect(\Delta^q_t)\cap(C_1e^{-C_2t},C_3t)=\emptyset$, and
\item
the number of the eigenvalues of $\Delta^q_t$
in the interval
$[0,C_1e^{-C_2t}]$ counted with their multiplicity is
equal to the number of zeros of $\omega$ of index $q$.
\endroster
\endproclaim

The above theorem states the existence of a gap in the spectrum of
$\Delta^q_t$, namely the open interval $(C_1e^{-C_2t},C_3t)$,
which widens to $(0,\infty)$ when $t\to\infty$.

Clearly $C_1,C_2,C_3$ and $T_0$ determine a constant $T\geq T_0$, so
that
for $t\geq T$, $1\in(C_1e^{-C_2t},C_3t)$
and therefore
$$
\spect(\Delta^q_t)\cap[0,C_1e^{-C_2t}]=\spect(\Delta^q_t)\cap[0,1]
$$
and
$$
\spect(\Delta^q_t)\cap[C_3t,\infty)=\spect(\Delta^q_t)\cap[1,\infty).
$$
For $t\geq T$ we denote by $\Omega^q_\tsm(M)$ the
finite dimensional subspace, generated by the eigenforms
of $\Delta^q_t$ corresponding to the eigenvalues of $\Delta^q_t$
smaller
than $1$.
The elliptic theory implies that these eigenvectors, a priori
elements
in the
$L^2$-completion of $\Omega^q(M)$, are actually in
$\Omega^q(M)$. Note
that
$d^q_t:\Omega^q_\tsm(M)\to\Omega^{q+1}_\tsm(M)$, so that
$\bigl(\Omega^*_\tsm(M),d^*_t\bigr)$ is a finite dimensional
cochain subcomplex
of $\bigl(\Omega^*(M),d^*_t\bigr)$. Clearly the
orthogonal complement
is also a closed subcomplex, we will denote by
$\bigl(\Omega^*_\tla(M),d^*_t\bigr)$.
One has the following orthogonal
decomposition
$$
\bigl(\Omega^*(M),d^*_t\bigr)
=\bigl(\Omega^*_\tsm(M),d^*_t\bigr)
\oplus
\bigl(\Omega^*_\tla(M),d^*_t\bigr),
$$
and $(\Omega^*_\tla(M),d^*_t)$ is acyclic.

Let $(\omega,g)$ be a Morse-Smale pair.
Recall that for each critical
point
$x\in\Cr(\omega)$ we have $\delta_x\in\maps(\Cr(\omega),\R)$ which takes the
value $1$ on $x$ and $0$ on all other critical points. Clearly
$\{\delta_x\mid x\in\Cr(\omega)\}$ is a base of the vector space
$\maps(\Cr(\omega),\R)$.
We equip $\maps(\Cr(\omega),\R)$ with the
unique scalar product which makes this base orthonormal.

The next result is an extension of Helffer-Sj\"ostrand theorem as formulated
in \cite{BFKM96},
but for closed one forms instead of functions.

\proclaim{Theorem 4}
Suppose $(\omega,g)$ is a Morse-Smale pair with $\rho(\omega,g)< \infty$ and $o$ are orientations as above.
Then there exists $T\geq 0$, depending
on $(\omega,g)$ so that for $t\geq T$
$$
\Int_t:\bigl(\Omega^*_\tsm(M),d^*_t\bigr)\to
\bigl(\maps(\Cr_*(\omega),\R),\partial^*_t\bigr)
$$
is an
isomorphism of cochain
complexes.
Moreover, there exists a family of
isometries
$R^q_t:\maps(\Cr_q(\omega),\R)\to\Omega^q_\tsm(M)$ of finite
dimensional
vector
spaces so that
$\Int_t\circ R^q_t=\id+O(1/t)$.
\endproclaim

The proof of Theorems~3 and 4 is similar to the one given in
\cite{BFKM96} or \cite{BFK} for Witten and Helffer-Sj\"ostrand theorems.
However for the
readers convenience we sketch the arguments in section~6.

For $t\geq T$ consider $E_{t,x}:=(\Int_t)^{-1}(\delta_x)\in\Omega^*_\tsm(M)$.
Clearly these forms provide a base for $\Omega^*_\tsm(M)$, and the
functions $t\mapsto I_{q,t}(x,y)$ are the unique functions which
satisfy
the formula
$$
d^q_t(E_{t,y})=\sum_{x\in\Cr_{q+1}(\omega)}I_{q+1,t}(x,y)E_{t,x}.
$$
Consequently the numbers $I_q(\tilde x,\tilde y)$ can be recovered from the
family
$(\Omega^*_\tsm(M),d^*_t)$ and the base
$\{E_{t,x}\mid x\in\Cr(\omega)\}$
by using the theory of Dirichlet series.

\head 2. The proof of Proposition 2\endhead

We will begin with few notations. Let $\bar\R:=\R\cup\{\pm\infty\}$,
equipped with the structure of a manifold with boundary via the diffeomorphism
$$
\bar\R\to[-1,1],\quad t\mapsto t(1+t^2)^{-\frac{1}{2}}.
$$
Choose a tubular neighborhood $V\subseteq TM$ of the zero section, such that
$$
(\exp,p):TM\supseteq V\to M\times M
$$
becomes a diffeomorphism onto its image. Here $\exp:TM\to M$ is defined
with respect to a Riemannian metric $g_0$ fixed once and for all.
Departing from the notation used in section~1 and in order to
remain as close as possible to the reference \cite{Sch93} we will write
$\gamma$ for an element of $\cixy:=\{\gamma\in C^\infty(\bar\R,M)\mid
\gamma(-\infty)=x,\gamma(\infty)=y\}$.

For $\gamma\in\cixy$ we have a well defined Sobolev space $\hgtm$,
cf \cite{Sch93} page~24 and a Sobolev embedding $\hgtm\subseteq
C^0(\gamma^*TM)$. Here $\gamma^*TM$ denotes the pull back of $TM\to M$
by $\gamma:\bar\R\to M$. So
$$
\hgu
:=\bigl\{\sigma\in\hgtm\bigm|\forall t\in\R:\sigma(t)\in V\bigr\}
$$
is an open neighborhood of $0\in\hgtm$. We set
$$
\varphi_\gamma:\hgu\to\czxy,
\quad\varphi_\gamma(\sigma):=\exp\circ\sigma,
$$
define
$$
\pxy:=\bigcup_{\gamma\in\cixy}\Img(\varphi_\gamma)
$$
and
$$
\exy:=\bigcup_{\sigma\in\pxy}\{\sigma\}\times\lstm
$$
and denote by $\pi:\exy\to\pxy$ the obvious projection. For
$\gamma\in\cixy$ let
$$
\psi_\gamma:\hgu\times\lgtm\to\exy,
\quad
\psi_\gamma(\sigma,\xi):=\big(\varphi_\gamma(\sigma),\Pt_1\xi\big),
$$
where $\Pt_1$ denotes the time 1 parallel transport along the geodesics
$s\mapsto\varphi_\gamma(s\sigma(t))$ with respect to the metric $g_0$.
The following facts are not hard to verify, cf Proposition~2.7 and
Proposition~2.9 in \cite{Sch93}.

\remark{Fact 1}
The maps $\varphi_\gamma$ \resp $\psi_\gamma$, $\gamma\in\cixy$, define an
atlas which provides a structure of smooth Hilbert manifold on
$\pxy$ \resp $\exy$ as well as a structure of a smooth Hilbert
vector bundle for $\pi:\exy\to\pxy$. These structures are, up to an
isomorphism, independent on the metric $g_0$. The inclusions
$\cixy\subseteq\pxy\subseteq\czxy$ are continuous maps, have dense images
and are homotopy equivalences.
\endremark

\remark{Fact 2}
Let $\omega\in\Cal Z^1(M)$ and $x,y\in\Cr(\omega)$. The map
$F:\pxy\to\exy$, defined by
$c\mapsto\bigl(c,\ddt c+(\grad_g\omega)\circ c\bigr)$
is a smooth section. If $x$ and $y$ are non-degenerate, then the zeros of $F$
are precisely the smooth mappings $\gamma:\R\to M$, satisfying
$$
\gamma'(t)=-(\grad_g\omega)(\gamma(t)),
\quad
\lim_{t\to-\infty}\gamma(t)=x
\quad\text{and}\quad
\lim_{t\to\infty}\gamma(t)=y.
$$
We will write $\Cal M(x,y):=F^{-1}(0)$, which identifies with
$W^-_x\cap W^+_y$. If $x,y$ are non-degenerated, then
$$
\pr_2\circ\psi_\gamma^{-1}\circ F\circ\varphi_\gamma:\hgu\to\lgtm
$$
is a Fredholm mapping of index $\index(x)-\index(y)$.
\endremark

\remark{Fact 3}
For $\gamma\in\Cal M(x,y)$ the differential of $F$ can be calculated
using the charts $\varphi_\gamma$ and $\psi_{\gamma}$. Precisely the
differential of $\psi_\gamma^{-1}\circ F\circ\varphi_\gamma$ at
$0$ is the linear map
$$
\aligned
T_0(\psi_\gamma^{-1}\circ F\circ\varphi_\gamma):\hgtm&\to\hgtm\times\lgtm
\\
T_0(\psi_\gamma^{-1}\circ F\circ\varphi_\gamma)(\xi)&=
\big(\xi,\nabla_{\dt}\xi+\nabla_\xi\grad_g\omega\big),
\endaligned
$$
where $\nabla$ denotes the connection on $\gamma^*TM$
induced from the Levi-Civita connection on $M$ provided by the
Riemannian metric $g_0$, and $\nabla_{\partial t}$ the induced connection
on $\gamma^*TM$.
\endremark

Suppose $\omega$ is a Morse form, $x\in\Cr(\omega)$ and let
$U_x:=B(x,\varepsilon)\setminus\overline{B(x,\frac\varepsilon2)}$. We set
$$
\Cal S_{U_x}:=\bigl\{A\in C^\infty(\es(TM))\bigm|
\supp(A-\id)\subseteq U_x\bigr\},
$$
where $\es(TM)$ denotes the endomorphisms of $TM$, which are symmetric with
respect to $g$. For a sequence $(\lambda_k)_{k\in\N}$ of real numbers and
$A\in C^\infty(\es(TM))$ we set
$$
\|A\|:=\sum_{k\geq0}\lambda_k
\sup|\underbrace{\nabla\cdots\nabla}_{\text{$k$-times}}A|,
$$
where $\nabla$ denotes the covariant differentiation induced from $g_0$, and
define $S_x:=\{A\in\Cal S_{U_x}\mid \|A\|<\infty\}$. One can choose
$(\lambda_k)_{k\in\N}$, such that $(S_x,\|\cdot\|)$ becomes an affine Banach
space, which is dense in $\Cal S_{U_x}$ with respect to the $L^2$-topology,
cf \cite{Sch93} and Lemma~5.1 in \cite{F88}. Finally let
$$
\Cal S^+_{U_x}:=
\bigl\{A\in\Cal S_{U_x}\bigm|\text{$A$ is pos\. def\. w\. r\. to $g$}\bigr\}
$$
and $S_x^+:=S_x\cap\Cal S^+_{U_x}$. Note that $S_x^+$ is an open neighborhood
of $\id$ in $S_x$. Consider the smooth mapping
$$
\Phi:S^+_x\times\pxy\to\exy,\quad
\Phi(A,c)=\bigl(c,\ddt c+(A\grad_g\omega)\circ c\bigr).
$$
Notice that $g_A(\cdot,\cdot):=g(A^{-1}\cdot,\cdot)$ is a Riemannian
metric, and $\Phi_A(\cdot):=\Phi(A,\cdot)=F_{g_A}$, since
$\grad_{g_A}\omega=A\grad_g\omega$. So for any $y\in\Cr(\omega)$ and any
$\gamma\in\cixy$
$$
\pr_2\circ\psi_\gamma^{-1}\circ\Phi_A\circ\varphi_\gamma:\hgu\to\lgtm
$$
is a Fredholm mapping of index $\index(x)-\index(y)$. If
$\Phi(A,\gamma)=0$ the differential at $(A,0)$ of the vertical part of
$\Phi$ in the chart given by $\varphi_\gamma$ and $\psi_\gamma$ is
$$
\align
D:S_x\times\hgtm&\to\lgtm\\
D(B,\xi)
&:=T_{(A,0)}\big(\pr_2\circ\psi_\gamma^{-1}\circ\Phi
\circ(\id,\varphi_\gamma)\big)(B,\xi)
\tag 2.1
\\&=
\nabla_{\dt}\xi+\nabla_\xi(A\grad_g\omega)+(B\grad_g\omega)\circ\gamma.
\endalign
$$

\proclaim{Lemma 1}
Let $\omega$ be a Morse form and $x,y\in\Cr(\omega)$.
Then $\Phi:S_x^+\times\pxy\to\exy$ intersects the zero section transversally.
\endproclaim

\demo{Proof}
Suppose $\Phi(A,\gamma)=0$. We have to show, that $D:S_x\times\hgtm\to\lgtm$
is onto. Since
$\pr_2\circ\psi_\gamma^{-1}\circ\Phi_A\circ\varphi_\gamma$ is a
Fredholm mapping, we see that $\Img(D)\subseteq\lgtm$
is a closed subspace of finite codimension. Suppose there exists
$0\neq\eta\in\lgtm$, such that
$$
\big\langle D(B,\xi),\eta\big\rangle=0,\quad
\forall(B,\xi)\in S_x\times\hgtm,
\tag 2.2
$$
where $\langle\cdot,\cdot\rangle$ denotes the inner product of
$\lgtm$. From \thetag{2.1} one gets
$$
\big\langle\nabla_{\dt}\xi+\nabla_\xi(A\grad_g\omega),\eta\big\rangle=0,
\quad\forall\xi\in\hgtm.
$$
The adjoint of $\xi\mapsto\nabla_{\dt}\xi+\nabla_\xi(A\grad_g\omega)$
is of the form $\eta\mapsto-\nabla_{\dt}\eta+K\eta$, for some
$K\in C^\infty(\end(\gamma^*TM))$.
So $\eta(t)\neq 0$ for all $t\in\R$ by the uniqueness result for ODEs.
Choose $t_0\in\R$, such that $\gamma(t_0)\in U_x$. From \thetag{2.2} and
\thetag{2.1} we also get
$$
\big\langle (B\grad_g\omega)\circ\gamma,\eta\big\rangle=0,
\quad\forall B\in S_x,
$$
and hence
$$
g\big((B_0\grad_g\omega)(\gamma(t_0)),\eta(t_0)\big)=0,
\quad\forall B_0\in\es(T_{\gamma(t_0)}M).
$$
But since $(\grad_g\omega)(\gamma(t_0))\neq 0$, we must have
$\eta(t_0)=0$, a contradiction. This verifies the surjectivity of $D$.
\qed
\enddemo

\proclaim{Lemma 2}
Let $x,y\in\Cr(\omega)$ and suppose $F:\pxy\to\exy$ is transversal to the
zero section. Then $i_x^-:\R^{\index(x)}\to W^-_x\subseteq M$ and
$i^+_y:\R^{n-\index(y)}\to W^+_y\subseteq M$ are transversal.
\endproclaim

\demo{Proof}
Suppose conversely, that they are not transversal at some point. Since every
point in the intersection  $W^-_x\cap W^+_y$ lies on a trajectory
$\gamma\in M(x,y)$, we may assume that this point is $\gamma(0)$.
Choose $\eta_0$ in the orthogonal complement of
$T_{\gamma(0)}W^-_x+T_{\gamma(0)}W^+_y$ and let $\eta$ be the unique vector
field along $\gamma$ satisfying
$$
-\nabla_{\dt}\eta+K\eta=0
\quad\text{and}\quad
\eta(0)=\eta_0,
$$
where $K$ is the pointwise adjoint of $\xi\mapsto\nabla_\xi\grad_g\omega$.
If we can show, that $\eta\in\lgtm$, then we would get a contradiction,
because
$$
\big\langle\nabla_{\dt}\xi+\nabla_\xi\grad_g\omega,\eta\big\rangle
=\big\langle\xi,-\nabla_{\dt}\eta+K\eta\big\rangle=0,
$$
but $\xi\mapsto\nabla_{\dt}\xi+\nabla_\xi\grad_g\omega$ is onto,
since $F$ is transversal to the zero section.

For $\xi_0\in T_{\gamma(0)}W^-_x$ let $\xi$ be the unique vector field
along $\gamma$ satisfying
$$
\nabla_{\dt}\xi+\nabla_\xi\grad_g\omega=0
\quad\text{and}\quad
\xi(0)=\xi_0.
$$
From
$\nabla_{\dt}\xi+\nabla_\xi\grad_g\omega=[\xi,\grad_g\omega]\circ\gamma$
it follows, that $\xi(t)\in T_{\gamma(t)}W^-_x$ for all $t$. So
$$
\ddt g(\xi,\eta)=g(\nabla_{\dt}\xi,\eta)+g(\xi,\nabla_{\dt}\eta)
=g(\nabla_{\dt}\xi+\nabla_\xi\grad_g\omega,\eta)=0,
$$
and hence $\eta(t)$ is orthogonal to $T_{\gamma(t)}W^-_x$ for all $t$.
A similar argument shows, that $\eta(t)$ is orthogonal to
$T_{\gamma(t)}W^+_y$ for all $t$, too.

Now consider the function $\alpha(t):=\frac12g(\eta(t),\eta(t))>0$.
Then
$$
\alpha'(t)
=g(K\eta(t),\eta(t))
=g\bigl(\eta(t),\nabla_{\eta(t)}\grad_g\omega\bigr)
=(\nabla_{\eta(t)}\omega)(\eta(t)).
$$
Since the Hessian of $\omega$ at $y$, see \thetag{1.3}, is negativ definite
on the orthogonal complement of stable manifold we find a constant $k>0$,
such that
$$
\alpha'(t)
=(\nabla_{\eta(t)}\omega)(\eta(t))
\leq-\tfrac12kg(\eta(t),\eta(t))
=-k\alpha(t),
$$
for large $t$. So $\ddt\ln\alpha(t)\leq-k$, hence
$\ln\alpha(t)\leq\ln(\alpha(0))-kt$ and
finally $\alpha(t)\leq\alpha(0)e^{-kt}$, for large $t$. So
we see that $\eta(t)$ converges exponentially to $0$ as $t\to\infty$.
A similar argument shows the exponential convergence for $t\to-\infty$.
This shows $\eta\in\lgtm$, and the proof is complete.
\qed
\enddemo

We are now in the position to give the

\demo{Proof of Proposition 2}
Recall first that a residual set in a complete metric space is a countable
intersection of open and dense sets. By Baire category theorem it is a
dense subset. Clearly a finite  intersection of residual sets is residual.

Next note, that every Riemannian metric on $M$ is of the form $g_A$ for a
unique positive definite $A\in C^\infty(\es(TM))$. We set
$G_x:=\{g_A\mid A\in S_x^+\}$ and
$$
G:=\prod_{x\in\Cr(\omega)}G_x
\subseteq\prod_{x\in\Cr(\omega)}\Cal G_{g,U_x}=\Cal G_{g,U}.
$$
Since $S^+_x$ is a Banach manifold which is dense in $\Cal S^+_{U_x}$ with
respect to the $L^2$-topology, the same is true fro $G\subseteq\Cal G_{g,U}$.
From Lemma~1 and Sard's theorem for Fredholm maps between
Banach manifolds, cf Proposition~2.24 in \cite{Sch93} it follows,
that for every $y\in\Cr(\omega)$ there exists a residual subset
$S_{x,y}'\subseteq S_x^+$, such that for any
$A\in S'_{x,y}$ the section $F_{g_A}:\pxy\to\exy$
intersects the zero section transversally. So
$$
G_x':=\bigl\{g_A\bigm| A\in\bigcap_{y\in\Cr(\omega)}S_{x,y}'\bigr\}
$$
is a residual subset of $G_x$, and Lemma~2 implies, that for any
$g'\in G_x'$ and any $y\in\Cr(\omega)$ the mappings $i^-_x$ and $i^+_y$
are transversal. So
$$
G':=\prod_{x\in\Cr(\omega)}G'_x\subseteq
\prod_{x\in\Cr(\omega)}G_x=G
$$
satisfies the statement of Proposition~2.
\qed
\enddemo

\head 3. The proof of Proposition 3\endhead

\proclaim{Lemma 3}
Let $(\omega,g)$ be a Morse-Smale pair. Then there exists a constant $C>0$,
such that
\roster
\item
For all $\tilde x,\tilde y\in\Cr(h)$ for which
$\Cal T(\tilde x,\tilde y)\neq\emptyset$ one has
$d(\tilde x,\tilde y)\leq C(h(\tilde x)-h(\tilde y))$.
\item
For all $\tilde x\in\Cr(h)$ and all $\tilde z\in W^-_{\tilde x}$ one has
$d(\tilde x,\tilde z)\leq\max\{C(h(\tilde x)-h(\tilde z)),1\}$.
\endroster
Here $d$ denotes the distance in $\tilde M$ given by the Riemannian
metric.
\footnote{
Actually the proof shows, that there exists a small ball
$B(\tilde x,\varepsilon)$, such that for
$\tilde z\in W_{\tilde x}^-\setminus B(\tilde x,\varepsilon)$ one has
$d_{W_{\tilde x}^-}(\tilde x,\tilde z)\leq
C\bigl(h(\tilde x)-h(\tilde z)\bigr)$, where $d_{W_{\tilde x}^-}$ denotes
the distance given by the induced Riemannian metric on $W_{\tilde x}^-$.
The only extra argument is needed in \thetag{3.4}, where one has to use the
fact that for every $\tilde y\in\Cr(h)$ every trajectory in
$B(\tilde y,\varepsilon)$ has length at most $2\varepsilon$.
}
\endproclaim

\demo{Proof}
For $r>0$ denote
$$
\tilde U_r:=\bigcup_{z\in\Cr(\omega)}\pi^{-1}\big(B(z,r)\big),
$$
where $B(z,r)$ denotes the open ball of radius $r$. Now choose
$\frac12\geq\varepsilon>0$, such that $\tilde U_\varepsilon$ is a disjoint
union of balls. Choose $C$, such that
$$
\frac4{\|(\grad_gh)(z)\|}\leq C,\quad
\text{for all $z\in\tilde M\setminus\tilde U_{\frac\varepsilon2}$.}
\tag 3.1
$$
Let $\tilde x,\tilde y\in\Cr(h)$ and $\gamma\in\Cal T(\tilde x,\tilde y)$,
parameterized by the value of $h$, cf Observation~5. So
$\gamma:[a,b]\to\tilde M$, where $a=h(\tilde y)$ and $b=h(\tilde x)$.

Suppose we have $[s,t]\subseteq[a,b]$. If
$\gamma([s,t])\subseteq\tilde M\setminus\tilde U_{\frac\varepsilon2}$ then
in view of \thetag{1.4} and \thetag{3.1} we get
$$
d\big(\gamma(s),\gamma(t)\big)
\leq\int_s^t|\gamma'(\sigma)|d\sigma
\leq\frac C4(t-s).
\tag 3.2
$$
If
$\gamma([s,t])\cap\partial\tilde U_{\frac\varepsilon2}\neq\emptyset$
and
$\gamma([s,t])\cap\partial\tilde U_\varepsilon\neq\emptyset$ then there
exists $s',t'\in[s,t]$, such that
$\gamma([s',t'])\subseteq\tilde M\setminus\tilde U_{\frac\varepsilon2}$,
$\gamma(s')\in\partial\tilde U_{\frac\varepsilon2}$ and
$\gamma(t')\in\partial\tilde U_\varepsilon$. So \thetag{3.2} yields
$$
\frac\varepsilon2
=d\big(\partial\tilde U_{\frac\varepsilon2},\partial\tilde U_\varepsilon\big)
\leq d\big(\gamma(s'),\gamma(t')\big)
\leq\frac C4|t'-s'|
\leq\frac C4|t-s|.
\tag 3.3
$$
This implies that there exist $a=s_0<t_0<s_1<t_1<\cdots<s_k<t_k=b$, such
that
$\gamma([t_i,s_{i+1}])\subseteq\tilde M\setminus\tilde U_{\frac\varepsilon2}$,
$\gamma((s_i,t_i))\subseteq\tilde U_\varepsilon$,
$\gamma([s_i,t_i])\cap\partial\tilde U_{\frac\varepsilon2}\neq\emptyset$
and $\gamma([s_i,t_i])\cap\partial\tilde U_\varepsilon\neq\emptyset$.
So \thetag{3.2} and \thetag{3.3} imply
$$
d\big(\gamma(t_i),\gamma(s_{i+1})\big)\leq C(s_{i+1}-t_i),
\quad
d\big(\gamma(s_i),\gamma(t_i)\big)\leq 2\varepsilon\leq C(t_i-s_i).
\tag 3.4
$$
Adding all these estimates together gives
$$
d(\tilde x,\tilde y)
\leq C(t_k-s_0)=C(b-a)=C\big(h(\tilde x)-h(\tilde y)\big).
$$
This proves part~\therosteritem1. To see part~\therosteritem2 notice, that
if $\tilde z$ does not lie in the component of $\tilde U_\varepsilon$
containing $\tilde x$, the argument above works and one gets
$d(\tilde x,\tilde z)\leq C(h(\tilde x)-h(\tilde z))$. If both lie in the
same component, one certainly has
$d(\tilde x,\tilde z)\leq2\varepsilon\leq 1$.
\qed
\enddemo

\proclaim{Corollary 1}
Let $(\omega,g)$ be a Morse-Smale pair. Then the following holds:
\roster
\item
For all $\tilde x,\tilde y\in\Cr(h)$ and $R\in\R$ the set
$\{\gamma\in\Gamma\mid
\Cal T(\gamma\tilde x,\tilde y)\neq\emptyset,[\omega](\gamma)\leq R\}$
is finite.
\item
Given $\tilde x,\tilde y\in\Cr(h)$, there exist only finitely many
$\tilde y_1,\dotsc,\tilde y_{k-1}\in\Cr(h)$, such that
$\Cal T(\tilde x,\tilde y_1)
\times\cdots\times\Cal T(\tilde y_{k-1},\tilde y)\neq\emptyset$.
\endroster
\endproclaim

\demo{Proof}
By Lemma~3\therosteritem1 and \thetag{1.1} the set $\{\gamma\tilde x\mid
\Cal T(\gamma\tilde x,\tilde y)\neq\emptyset, [\omega](\gamma)\leq R\}
\subseteq\tilde M$ is bounded, and since it is discrete too, it must be
finite, for $\tilde M$ is a complete Riemannian manifold.
Statement~\therosteritem1 follows immediately. Part~\therosteritem2
follows from a similar argument, $\index(\tilde y_i)>\index(\tilde y_{i+1})$
and $h(\tilde x)\geq h(\tilde y_i)\geq h(\tilde y)$.
\qed
\enddemo

\demo{Proof of Proposition 3}
We will only prove part \therosteritem1, the proof of \therosteritem2 is
similar. First we will show that $\Cal B(\tilde x,\tilde y)$ is closed.
For notational simplicity set $a:=h(\tilde y)$ and $b:=h(\tilde x)$.
Suppose $\gamma_n\in\Cal B(\tilde x,\tilde y)$ converge uniformly to
$\gamma_\infty\in C^0\big([a,b],\tilde M\big)$. Clearly
the conditions in Observation~5\therosteritem1 hold for $\gamma_\infty$, too.
Since $\gamma_\infty([a,b])$ is compact and $\Cr(h)$ is discrete and because
of $h(\gamma_\infty(s))=a+b-s$, there are only finitely many
$s_i\in[a,b]$, with $\gamma_\infty(s_i)\in\Cr(h)$. If
$\gamma_\infty(s)\notin\Cr(h)$ then, for large $n$, the same holds for
$\gamma_n(s)$ and \thetag{1.4} follows. So $\Cal B(\tilde x,\tilde y)$
is closed.

Lemma~3\therosteritem2 implies, that
$\{\gamma([a,b])\mid\gamma\in\Cal B(\tilde x,\tilde y)\}\subseteq
\tilde M$ is bounded and since $\tilde M$ is a complete Riemannian
manifold, its closure is compact. In view of the
theorem of Arzela-Ascoli it remains to show that
$\Cal B(\tilde x,\tilde y)$ is equicontinuous. So let $\varepsilon>0$ small
and let $C_\varepsilon$ denote the constant $C$ we have constructed in the
proof of Lemma~3, which actually depended on $\varepsilon$. Set
$\delta:=\frac\varepsilon{C_\varepsilon}$ and suppose $s_0,s_1\in[a,b]$,
with $|s_1-s_0|\leq\delta$. We have to show
$$
d\big(\gamma(s_0),\gamma(s_1)\big)\leq2\varepsilon,
\quad\text{for all $\gamma\in\Cal B(\tilde x,\tilde y)$.}
\tag 3.5
$$
If $\gamma([s_0,s_1])\subseteq\tilde M\setminus\tilde U_{\frac\varepsilon2}$
this follows from \thetag{3.2}. If
$\gamma([s_0,s_1])\cap\tilde U_{\frac\varepsilon2}\neq\emptyset$, we must
have $\gamma([s_0,s_1])\subseteq\tilde U_\varepsilon$, for otherwise we
get a contradiction to \thetag{3.3}. But since the diameter of each component
of $\tilde U_\varepsilon$ is $2\varepsilon$, \thetag{3.5} follows in this
case, too.
\qed
\enddemo

\head 4. The proof of Theorem 1\endhead

For didactical reasons the proof will be given first in the particular
case that the set of all critical values, $h(\Cr(h))$ is a discrete subset
of $\R$, \ie $\omega$ has degree of rationality 1. Then we will show, that
the same arguments properly modified hold in the general case as well.
For the case where $\omega$ is exact the proof below is similar to the one
in \cite{BFK}.

\subhead 4.1 Some notations\endsubhead
Let  $\cdots>c_i>c_{i-1}>\cdots$, $i\in\Z$ denote the set of all critical
values of $h$. Choose $\epsilon_i>0$ small enough, so that
$c_i-\epsilon_i>c_{i-1}+\epsilon_{i-1}$, for all $i\in\Z$. Denote,
see Figure~1,
$$
\align
\Cr(i):=&\Cr(h)\cap h^{-1}(c_i),\\
M_i:= &h^{-1}(c_i),\\
M^\pm_i:=&h^{-1}(c_i\pm\epsilon_i)\quad\text{and}\\
M(i):=&h^{-1}(c_{i-1},c_{i+1}).
\endalign
$$
In view of Observation~2, $\Cr(i)$ is always a finite set, even when $\omega$
has degree of rationality greater than 1.

To keep the notation simpler we will denote the critical points of of
$h$ by $x,y,\dotsc$ instead of $\tilde x,\tilde y,\dotsc$.
There is no danger of confusion since the critical points of $\omega$
will not appear in this section. For any $x\in\Cr(i)$ denote,
see Figure~1,
$$
\align
S^\pm_x:=&W_x^\pm\cap M^\pm_i,\\
\bold S_x:=&S^+_x\times S^-_x,\\
W_x^\pm(i):=&W_x^\pm\cap M(i)\quad\text{and}\\
\bold SW_x(i):=&S_x^+\times W_x^-(i).
\endalign
$$
It will be convenient to write
$$
\align
S^\pm_i:=&\bigcup_{x\in\Cr(i)}S^\pm_x,\\
\bold S_i:=&\bigcup_{x\in\Cr(i)}S_x,\\
W^\pm(i):=&\bigcup_{x\in\Cr(i)}W_x^\pm(i)\quad\text{and}\\
\bold SW(i):=&\bigcup_{x\in\Cr(i)}\bold SW_x(i).
\endalign
$$

$$
\gather
\epsfxsize=4in\epsfbox{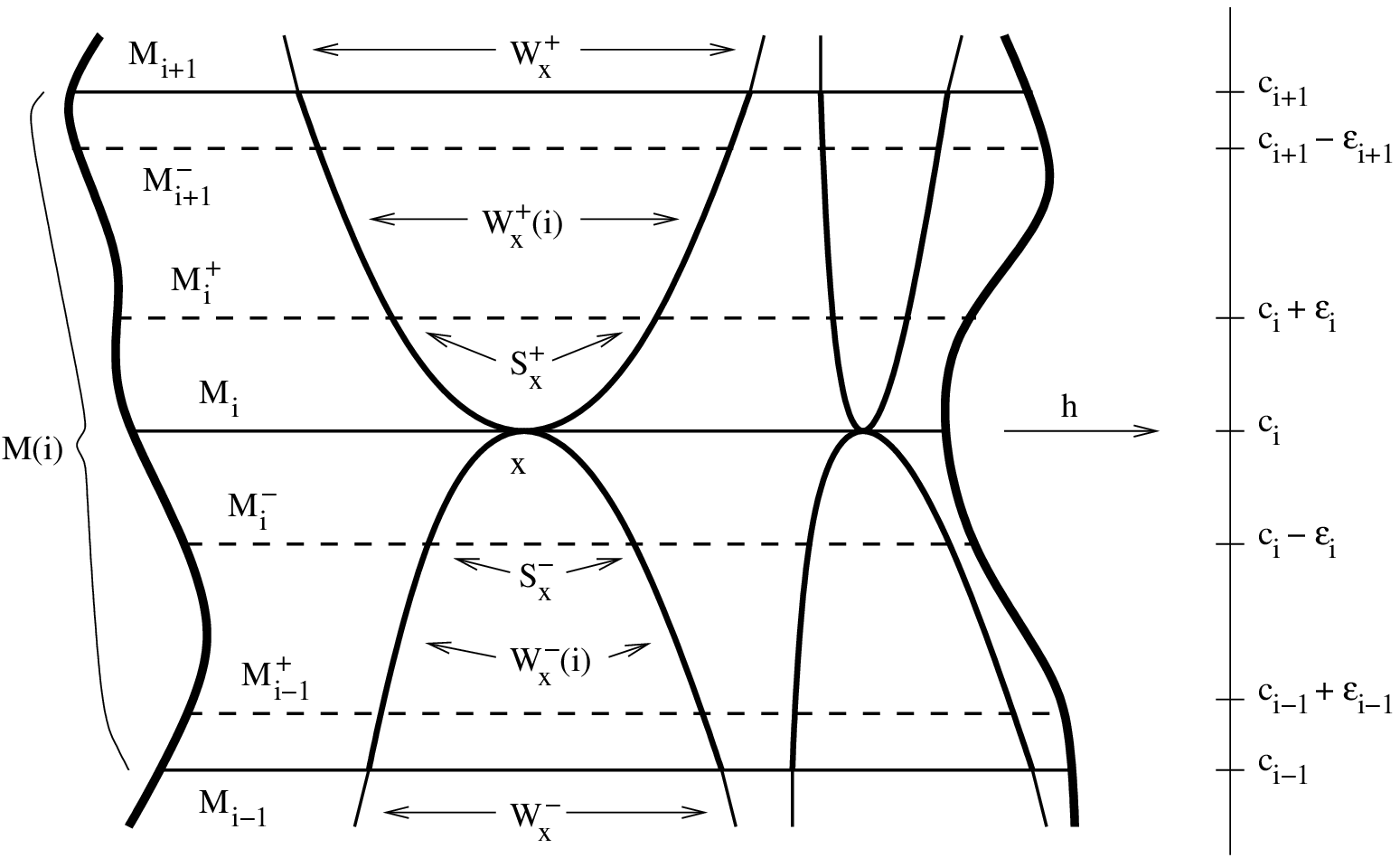}
\\
\text{Figure 1}
\endgather
$$

\remark{Observation 6}
We have:
\roster
\item
$\bold S_i\subseteq M^+_i\times M^-_i$.
\item
$\bold SW(i)\subseteq S^+_i\times W^-(i)\subseteq M^+_i\times M(i)$.
\item
$M^\pm_i$ are smooth manifolds of dimension $n-1$, where $n=\dim(\tilde M)$.
\item
$M(i)$ is a smooth manifold of dimension $n$, actually an open set in
$\tilde M$.
\item
$M_i$ is not a manifold, however
$\dot M_i:=M_i\setminus\Cr(i)$ and
$\dot M_i^\pm:=M_i^\pm\setminus S_i^\pm$
are are smooth manifolds of dimension $n$, actually submanifolds of
$\tilde M$.
\endroster
\endremark

Let $\Phi_t$ be the flow associated to the vector field
$-\grad_g h/||\grad_g h||^2$ on $\tilde M\setminus\Cr(h)$ and consider
the diffeomorphisms, see Figure~2,
$$
\psi_i:M^-_i\to M^+_{i-1},\quad
\psi_i(x):=\Phi_{c_i-c_{i-1}-\epsilon_i-\epsilon_{i-1}}(x)
$$
and
$$
\varphi_i^\pm:\dot M^\pm_i\to\dot M_i,\quad
\varphi_i^\pm(x):=\Phi_{\pm\epsilon_i}(x),
$$
as well as the submersion
$$
\varphi(i):M(i)\setminus\bigl(W^+(i)\cup W^-(i)\bigr)\to\dot M_i,\quad
\varphi(i)(x):=\Phi_{h(x)-c_i}(x).
$$

\remark{Observation 7} $\varphi_i^\pm$ and $\varphi(i)$ extend to
continuous maps
$$
\varphi_i^\pm:M^\pm_i\to M_i\quad\text{and}\quad\varphi(i):M(i)\to M_i.
$$
\endremark

$$
\gather
\epsfxsize=4in\epsfbox{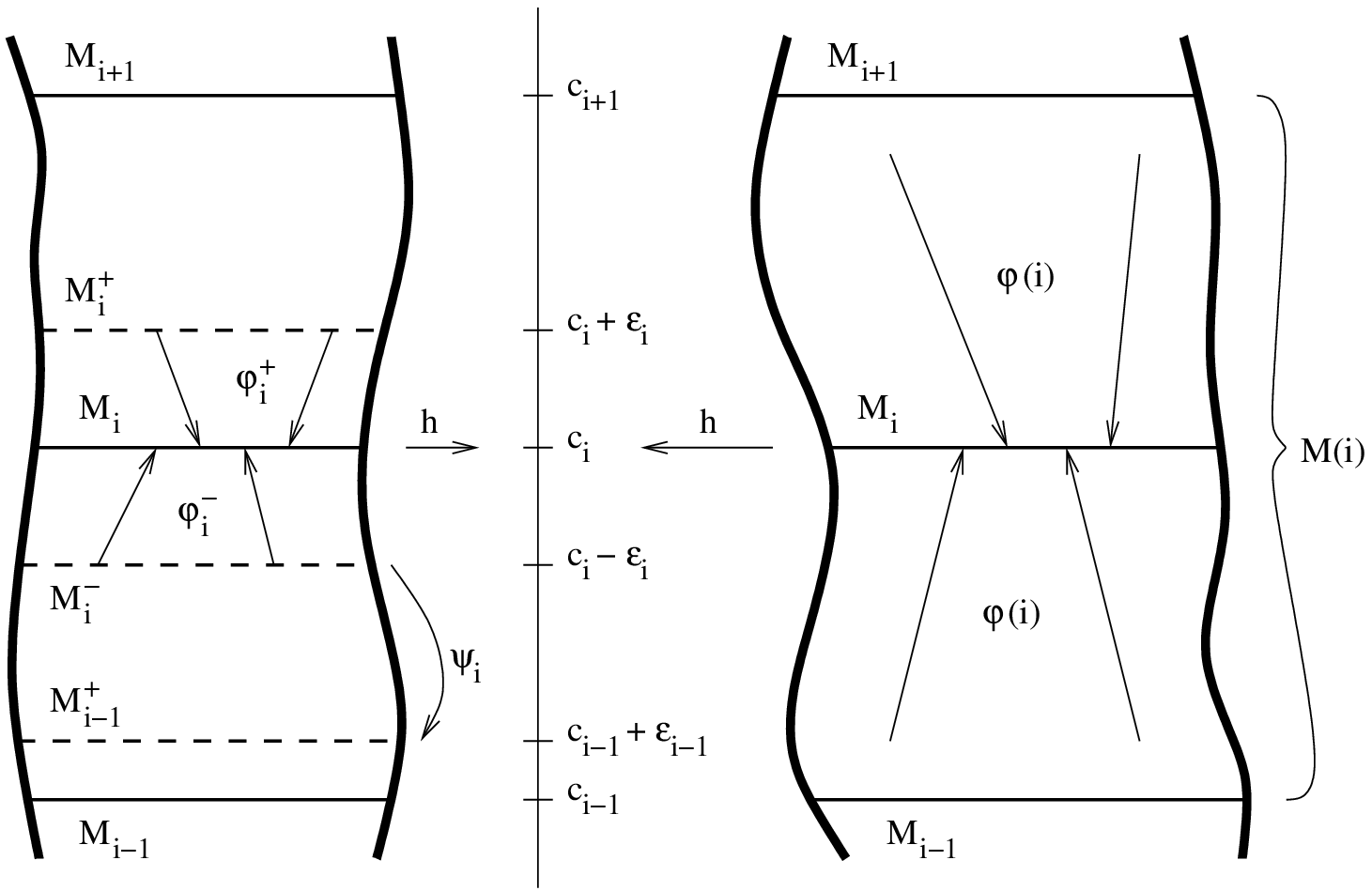}
\\
\text{Figure 2}
\endgather
$$

Define
$$
P_i:=\bigl\{(x,y)\in M_i^+\times M_i^-\bigm|
\varphi^+_i(x)=\varphi^-_i(y)\bigr\},
$$
and denote by $p^\pm_i:P_i\to M_i^\pm$ the canonical projections.
One can verify the following

\remark{Observation 8}
$P_i$ is a smooth $(n-1)$-dimensional manifold with boundary
(smooth submanifold of $M^+_i\times M^-_i$), whose boundary $\partial P_i$
is diffeomorphic to $\bold S_i\subset M^+_i\times M^-_i$.
Precisely we have
\roster\widestnumber\item{(P2)}
\item"(P1)"
$p^\pm_i:P_i\setminus\partial P_i\to\dot M_i^\pm$ are diffeomorphisms, and
\item"(P2)"
the restriction of $p^+_i\times p^-_i$ to $\partial P_i$
is a diffeomorphism onto $\bold S_i$, each $p^\pm_i$ restricted
to $\partial P_i$ identifies with the projection onto $S^\pm_i$.
\endroster
\endremark

Next we define
$$
Q(i):=\bigl\{(x,y)\in M_i^+\times M(i)\bigm|
\varphi^+_i(x)=\varphi(i)(y)\bigr\},
$$
or equivalently, $Q(i)$ consists of pairs of points $(x,y)$, $x\in M_i^+$,
$y\in M(i)$, which lie on the same (possibly broken) trajectory. Moreover
let $l_i:Q(i)\to M^+_i$ and $r_i:Q(i)\to M(i)$ denote the canonical
projections. One can verify the following

\remark{Observation 9}
$Q(i)$ is a smooth $n$-dimensional manifold with boundary (smooth
submanifold of $M^+_i\times M(i)$), whose boundary $\partial Q(i)$ is
diffeomorphic to $\bold SW(i)\subset M^+_i\times M(i)$. Precisely we have
\roster\widestnumber\item{(Q2)}
\item"(Q1)"
$l_i:Q(i)\setminus\partial Q(i)\to\dot M^+_i$ is a smooth bundle with
fiber an open segment, and
$r_i:Q(i)\setminus\partial Q(i)\to M(i)\setminus W^-(i)$ is a diffeomorphism.
\item"(Q2)"
The restriction of $l_i\times r_i$ to $\partial Q(i)$
is a diffeomorphism onto $\bold SW(i)$, \ie $l_i$ \resp $r_i$ restricted
to $\partial Q(i)$ identifies with the projection onto
$S^{+}_i$ \resp $W^-(i)$.
\endroster
\endremark

Since $P_i$ and $Q(i)$ are smooth manifolds with boundaries
$$
\Cal P_{r, r-k}:=P_r\times P_{r-1}\times\cdots\times P_{r-k}
$$
and
$$
\Cal P_r(r-k):=P_r\times\cdots\times P_{r-k+1}\times Q(r-k)
$$
are smooth manifolds with corners.

\subhead 4.2 The proof for degree of rationality 1\endsubhead
The proof of Theorem~1 will be based on the following recognition
method for a smooth manifold with corners.

\remark{Observation 10}
If $\Cal P$ is a smooth manifold with corners, $\Cal O,\Cal S$ smooth
manifolds, $p:\Cal P\to\Cal O$ and $s:\Cal S\to\Cal O$ smooth maps so that
$p$ and $s$ are transversal ($p$ is transversal to $s$ if its restriction to
each $k$-boundary $\Cal P_k$ is transversal to $s$), then $p^{-1}(s(\Cal S))$
is a smooth submanifold with corners of $\Cal P$.
\endremark

\demo{Proof of Theorem 1\therosteritem1}
We want to verify that $\Cal B(x,y)$ is a smooth manifold with corners. Let
$x\in\Cr(r+1)$ and $y\in\Cr(r-k-1)$, $k\geq -2$. If $k=-2$ the statement
is empty, if $k=-1$ there is nothing to check, so we suppose $k\geq 0$.

We consider $\Cal P=\Cal P_{r,r-k}$ as defined above,
$\Cal O:=\prod^{r-k}_{i=r}(M_i^+\times M_i^-)$ and
$\Cal S:=S^-_x\times M^-_r\times\cdots\times M^-_{r-k+1}\times S^+_y$.
In order to define the maps $p$ and $s$ we consider
$$
\omega_i:M^-_i\to M^-_i\times M^+_{i+1},\quad
\omega_i(x):=(x,\psi_i(x))
$$
and
$$
\tilde p_i:P_i\to M^+_i\times M^-_i,\quad
\tilde p_i(y):=\bigl(p^+_i(y),p^-_i(y)\bigr).
$$
We also denote by $\alpha:S^-_x\to M^+_r$ \resp $\beta:S^+_y\to M^-_{k-r}$ the
restriction of $\psi_{r+1}$ \resp $\psi_{r-k}^{-1}$ to $S_x^-$ \resp $S_y^+$.
Finally we set, see Diagram~1,
$$
s:=\alpha\times\omega_r\times\cdots\times\omega _{r-k+1}\times\beta:
\Cal S\to\Cal O
$$
and
$$
p:=\tilde p_r\times\cdots\times\tilde p_{r-k}:
\Cal P\to\Cal O.
$$

The verification of the transversality of $p$ and $s$ follows easily from
(P1), (P2) and the Morse-Smale condition, as we will explain in section~4.3
below. It is easy to see that $p^{-1}(s(\Cal S))$ identifies to
$\Cal B(x,y)$ as topological spaces and we leave this verification to the
reader. The compactness of $\Cal B(x,y)$ is stated in Proposition~3.
\qed
\enddemo

\demo{Proof of Theorem 1\therosteritem2}
Consider the set $X:=\hat W^-_x$, the map $\hat i_x:X=\hat W^-_x\to\tilde M$
and $\hat h:X\to\R$, see Definition~5. For any positive integer $k$, denote
by $X(k):=\hat i_x^{-1}(M(k))$. First we will topologize $X(k)$ and put on it
a structure of smooth manifold with corners, so that the restriction of
$\hat i_x$ and of $\hat h_x$ to $X(k)$ are smooth maps. Second we check
that $X(k)$ and $X(k')$ induce on the intersection $X(k)\cap X(k')$ the
same topology and the same smooth structure. These facts imply that $X$
has a canonical structure of smooth manifold with corners and that $\hat i_x$
is a smooth map. The properness of $\hat h_x$ follows from the compactness of
$\hat h^{-1}(s)$, which is in fact the space $\Cal B(x;s)$ whose compactness
is stated in Proposition~3.

To accomplish first step we proceed in exactly the same way as in the
proof of part~\therosteritem1. Suppose $x\in\Cr(r-1)$. Consider
$\Cal P:=\Cal P_r(r-k)$,
$\Cal O:=\prod_{i=r}^{r-k}(M_i^+\times M_i^-)$ and
$\Cal S:= S^-_x\times M^-_k\times\cdots\times M^-_{r-k+1}$.
Define, cf Diagram~2,
$$
p:=\tilde p_r\times\cdots\times\tilde p_{r-k+1}\times l_{r-k}:
\Cal P\to\Cal O
$$
and
$$
s:=\alpha\times\omega_r\times\cdots\times\omega_{r-k+1}:\Cal S\to\Cal O.
$$
The verification of the transversality follows from (P1), (P2), (Q1), (Q2)
and the Morse-Smale condition, as will be explained in section~4.3, below.
It is easy to see and left to the reader, that $p^{-1}(s(\Cal S))$
identifies to $X(r-k)$. The second step is more or less straightforward,
so it will be left again to the reader.
\qed
\enddemo

\subhead 4.3 The transversality of $p$ and $s$\endsubhead
Consider the diagrams:
$$
\gather
\xymatrix @C=-1pt{
S^-_{r+1}\ar[d]_\alpha
&&&M^-_r\ar[dl]_\id\ar[dr]^{\psi_r}
&&&&{M^-_{r-1}}\ar[dl]_\id
&\cdots
&{M^-_{r-k+1}}\ar[dr]^{\psi_{r-k+1}}
&&&S^+_{r-k-1}\ar[d]^\beta
\\
M^+_r
&&{M^-_{r}}
&&{M^+_{r-1}}
&&M^-_{r-1}
&&\cdots
&&{M^+_{r-k}}
&&{M^-_{r-k}}
\\
&P_r\ar[ul]^-{p^+_{r}}\ar[ur]_{p^-_{r}}
&&&&P_{r-1}\ar[ul]^{p^+_{r-1}}\ar[ur]_{p^-_{r-1}}
&&&\cdots
&&&P_{r-k}\ar[ul]^{p^+_{r-k}}\ar[ur]_{p^-_{r-k}}
}
\\
\text{Diagram 1}
\endgather
$$

$$
\gather
\xymatrix @C=3pt{
S^-_{r+1}\ar[d]_\alpha
&&&M^-_r\ar[dl]_\id\ar[dr]^{\psi_r}
&&&&{M^-_{r-1}}\ar[dl]_\id
&\cdots
&{M^-_{r-k+1}}\ar[dr]^{\psi_{r-k+1}}
\\
M^+_r
&&{M^-_{r}}
&&{M^+_{r-1}}
&&M^-_{r-1}
&&\cdots
&&{M^+_{r-k}}
\\
&P_r\ar[ul]^-{p^+_{r}}\ar[ur]_{p^-_{r}}
&&&&P_{r-1}\ar[ul]^{p^+_{r-1}}\ar[ur]_{p^-_{r-1}}
&&&\cdots
&&Q(r-k)\ar[u]_{l_{r-k}}
}
\\
\text{Diagram 2}
\endgather
$$

$$
\gather
\xymatrix @C=3pt{
&M^-_r\ar[dl]_{\id}\ar[dr]^{\psi_r}
&&&&{M^-_{r-1}}\ar[dl]_{\id}
&\cdots
&{M^-_{r-k+1}}\ar[dr]^{\psi_{r-k+1}}
&&&S^+_{r-k-1}\ar[d]^\beta
\\
{M^-_r}
&&{M^+_{r-1}}
&&M^-_{r-1}
&&\cdots
&&{M^+_{r-k}}
&&{M^-_{r-k}}
\\
S^-_r\ar[u]^i
&&&P_{r-1}\ar[ul]^{p^-_{r-1}}\ar[ur]_{p^+_{r-1}}
&&&\cdots
&&&P_{r-k}\ar[ul]^{p^+_{r-k}}\ar[ur]_{p^-_{r-k}}
}
\\
\text{Diagram 3}
\endgather
$$

$$
\gather
\xymatrix @C=3pt{
S^-_{r+1}\ar[d]_\alpha
&&&{M^-_r}\ar[dl]_{\id}\ar[dr]^{\psi_r}
&&&&{M^-_{r-1}}\ar[dl]_{\id}
&\cdots
&{M^-_{r-k+1}}\ar[dr]^{\psi_{r-k+1}}
\\
{M^+_r}
&&{M^-_{r}}
&&{M^+_{r-1}}
&&M^-_{r-1}
&&\cdots
& &{M^+_{r-k}}
\\
&{P_r}\ar[ul]^{p^+_r}\ar[ur]_{p^-_r}
&&&&P_{r-1}\ar[ul]^{p^+_{r-1}}\ar[ur]_{p^-_{r-1}}
&&&\cdots
&&S^+_{r-k}\ar[u]_i
}
\\
\text{Diagram 4}
\endgather
$$

$$
\gather
\xymatrix @C=3pt{
&{M^-_r}\ar[dl]_\id\ar[dr]^{\psi_r}
&&&&{M^-_{r-1}}\ar[dl]_\id
&\cdots
&{M^-_{r-k+1}}\ar[dr]^{\psi_{r-k+1}}
\\
{M^-_r}
&&{M^+_{r-1}}
&&M^-_{r-1}
&&\cdots
&&{M^+_{r-k}}
\\
S^-_r\ar[u]^i
&&&P_{r-1}\ar[ul]^{p^+_{r-1}}\ar[ur]_{p^-_{r-1}}
&&&\cdots
&&Q(r-k)\ar[u]_{l_{r-k}}
}
\\
\text{Diagram 5}
\endgather
$$
For each of these diagrams denote by $\Cal P$ \resp $\Cal O$ \resp
$\Cal S$ the product of the manifolds on the third  \resp second
\resp first row and let $p:\Cal P\to\Cal O$ \resp
$s:\Cal S\to\Cal O$ denote the product of the maps from the third to
the second row \resp from the first to the second row.
Clearly $\Cal P$ is a smooth manifold with corners. Denote by $\Cal P_0$
the interior of $\Cal P$, and by
$p_0:\Cal P_0\to\Cal O$ the restriction of $p$ to $\Cal P_0$.

We refer to the  statement \lq $p_0$ is transversal to $s$\rq\ with
$p_0$ and $s$ obtained from Diagram~$j$ as $T^j_{r,k}$, $j=1,\dotsc,5$.
Since all arrows but $\alpha$, $\beta$, $l_{r-k}$ and $i$
are open embeddings, the properties $T^2_{r, k}$ and $T^5_{r, k}$
follow. $T^1_{r,k}$ \resp $T^3_{r,k}$ \resp $T^4_{r,k}$ follow from
the transversality of $W^-_{r+1}$ and $W^+_{r-k-1}$ \resp
$W^-_r$ and $W^+_{r-k-1}$ \resp $W^-_{r+1}$ and $W^+_{r-k}$, \ie the
Morse-Smale condition.

Note that if $a_i:\Cal A_i\to\Cal B_i$ and $c_i:\Cal C_i\to\Cal B_i$
are transversal, $\Cal B_i$, $\Cal C_i$ smooth manifolds, $\Cal A_i$
smooth manifold with corners, $i=1,2$ then
$$
a_1\times a_2:\Cal A_1\times\Cal A_2\to\Cal B_1\times\Cal B_2
\quad\text{and}\quad
c_1\times c_2:\Cal C_1\times\Cal C_2\to\Cal B_1\times\Cal B_2
$$
are transversal, too. So in view of (P2) and (Q2)
it is easy to see, that the transversality of $p$ and $s$ obtained
from the diagram $1$ \resp $2$ can be derived
from the validity of the statements
$T^1_{r, k}$, $T^3_{r,k}$, $T^4_{r,k}$ \resp
$T^2_{r, k}$, $T^4_{r,k}$, $T^5_{r,k}$ for various $r,k$.

\remark{Observation 11}
If in Diagrams~1--5 above $\psi_r$  are only open embeddings
rather than diffeomorphisms, the transversality of $p$ and $s$ still
holds from the same reasons.
\endremark

\subhead 4.4 The general case\endsubhead
We start with the following

\definition{Definition 7 (Relevant critical points and values)}
\roster
\item
Let $x,y\in\Cr(h)$. Then $z\in\Cr(h)$ \resp the real number $h(z)$ is called
$(x,y)$-relevant critical point \resp $(x,y)$-relevant critical value,
if there exist $y_0,\dotsc,y_r\in\Cr(h)$ and $0\leq i_0\leq r$,
such that $y_0=x$, $y_r=y$, $y_{i_0}=z$ and such that
$$
\Cal T(y_0,y_1)\times\cdots\times\Cal T(y_{r-1},y_r)\neq\emptyset.
$$
\item
Let $x\in\Cr(h)$. Then $z\in\Cr(h)$ \resp the real number $h(z)$ is a called
$x$-relevant critical point \resp $x$-relevant critical value, if it is
$(x,y)$-relevant for some $y\in\Cr(h)$.
\endroster
\enddefinition

From Corollary~1 we immediately get the following

\remark{Observation 12}
Let $x,y\in\Cr(h)$. Then there are only finitely many $(x,y)$-relevant
critical points and values. Moreover the set of $x$-relevant critical values
is a discrete set of real numbers, bounded from above by $h(x)$. Note
that if $\omega$ has degree of rationality bigger than 1, the set of all
critical values is not discrete, but it still has measure 0, by Sard's
theorem.
\endremark

Let us consider $x\in\Cr(h)$ and denote by $h(x)=c_0>c_{-1}>\cdots$
the discrete set of $x$-relevant critical values. We choose
$\epsilon_i$ as above with the additional property that
$c_i\pm\epsilon_i$ are regular values. We proceed as in the previous
case but with care.
\roster
\item
In the definition of $S^\pm_i,\bold S_i,W^\pm(i)$ and $\bold SW(i)$ the
union should be taken only over critical points in
$\Cr_x(i):=\{y\in\Cr(i)\mid\text{$y$ is $x$-relevant}\}$.
\item
The diffeomorphisms $\psi_i$ and $\varphi_i^\pm$ are only partially
defined with maximal domains open sets in $M^-_i$ and
$\dot M^\pm_i$ but still diffeomorphisms onto their images, the submersion
$\varphi(i)$ with maximal domain an open set and the continuous extensions
$\varphi_i^\pm$ and $\varphi(i)$ partially defined with maximal domains open
sets.
\item
The sets $P_i$ will involve only pairs $(x,y)$ with $x$ in the
domain of $\varphi^-_i$ and $y$ in the domain of $\varphi^+_i$ and
$Q(i)$ will involve only pairs $(x,y)$, with $x$ in the domain of
$\varphi^+_i$ and $y$ in the domain of $\varphi(i)$. They remain however
manifolds with boundary.
\item
The conclusions (P2) and (Q2) remain the same and in (P1) \resp (Q1)
diffeomorphism \resp smooth bundle are replaced by open embedding \resp
submersion.
\endroster
With these specifications the proof is a word by word repetition of
the proof in the case of degree of rationality 1.

\head 5. The proof of Theorem 2\endhead

Let $\pi:\tilde M\to M$ be a covering corresponding to $[\omega]$.
Recall that  $\Cal X_q$ \resp $\Cr_q(\omega)$ denote the set of critical
points of $h$ \resp of $\omega$ of index $q$. $\Gamma$ acts freely on
$\Cal X_q$ with quotient set $\Cr_q(\omega)$. So
$C^q=\maps(\Cr_q(\omega),\C)$ can be identified via $\pi^*$ with the
$\Gamma$-invariant functions $\Cal X_q\to\C$. Moreover
$\pi^*:\Omega^q(M;\C)\to\Omega^q(\tilde M;\C)$ provides an identification
of $\Omega^q(M;\C)$ with the $\Gamma$-invariant $q$-forms on $\tilde M$.

The following proposition is a corollary of Theorem~1 and will be the main
tool in the proof of Theorem~2.

\proclaim{Proposition 4}
Let $s\in\C$ with $\re(s)>\rho(\omega,g)$. Then the following holds:
\roster
\item
For any $\alpha\in\Omega^q(M;\C)$ and any $\tilde x\in\Cal X_q$ the integral
$$
\Int_s(\alpha)(\tilde x):=\int_{W^-_{\tilde x}}e^{sh^{\tilde x}}\pi^*\alpha
\tag 5.1
$$
converges absolutely,
\footnote{
Recall that for an oriented $n$-dimensional manifold $N$ and 
$\alpha\in\Omega^n(N;\C)$ one has 
$|\alpha|:=|a|\Vol\in\Omega^n(M)$,
where $\Vol\in\Omega^n(N)$ is any volume form and 
$a\in C^\infty(N,\C)$ is the
unique function satisfying $\alpha=a\cdot\Vol$. 
The integral $\int_N\alpha$ is called absolutely convergent, if
$\int_N|\alpha|$ converges.
} 
does only depend on $\pi(\tilde x)$ and defines a
surjective linear map $\Int_s:\Omega^q(M;\C)\to C^q$.
\item
For any $\alpha\in\Omega^q(M;\C)$ and any $\tilde x\in\Cal X_{q+1}$ one has
$$
\Int_s\bigl(d^q_s(\alpha)\bigr)(\tilde x)
=\sum_{\tilde y\in\Cal X_q}
I_{q+1}(\tilde x,\tilde y)e^{-sH(\tilde x,\tilde y)}
\Int_s(\alpha)(\tilde y)
\tag 5.2
$$
\endroster
\endproclaim

\demo{Proof}
We start with part~\therosteritem1. Consider $F_q(M)\to M$ the smooth bundle
of orthonormal $q$-frames which is a compact smooth manifold. A
differential form $\alpha\in\Omega^q(M;\C)$ induces a smooth function
$\langle\alpha\rangle:F_q(M)\to\C$ which is bounded by a positive constant
$C_\alpha$, \ie $|\langle\alpha\rangle(\tau)|\leq C_\alpha$ for every
orthonormal frame $\tau$. Then we have
$$
\int_{W^-_{\tilde x}}\bigl|e^{sh^{\tilde x}}
\pi^*\alpha\bigr|\Vol_{W_{\tilde x}^-}
=\int_{\R^q}e^{\re(s)h^x}\bigl|(i^-_x)^*\alpha\bigr|\Vol_{g^x}
\leq C_\alpha\int_{\R^q}e^{\re(s)h^x}\Vol_{g^x},
$$
hence the convergence of the integral \thetag{5.1}
insured by the definition of $\rho(\omega,g)$.

To verify the surjectivity of $\Int_s$ we construct for each
$x\in\Cr_q(\omega)$ a smooth one parameter family of differential forms
$\alpha^x_\lambda\in\Omega^q(M;\C)$, $\lambda\in[0,\epsilon]$ with the
following properties:
\roster
\item
$\lim_{\lambda\to 0}\int_{W^-_{\tilde x}}e^{sh^{\tilde x}}
\pi^*(\alpha^x_\lambda)=1$, for any
$\tilde x\in\tilde M$ with $\pi(\tilde x)=x$.
\item
If $x'\neq x$ but $\index(x)=\index(x')$ then
$\lim_{\lambda\to 0}\int_{W^-_{\tilde x'}}e^{sh^{\tilde x'}}
\pi^*(\alpha^x_\lambda)=0$.
\endroster
It is then clear, that by taking $\lambda$ small enough
$\Int_s(\alpha^x_\lambda)$, $x\in\Cr_q(\omega)$ are linearly independent,
hence a base of $C^q$, and therefore $\Int_s$ is surjective.

Now let us  describe the construction of the family $\alpha^x_\lambda$.
We use coordinates $(t_1,\dotsc,t_r)$ to parameterize points in $\R^r$ and
denote by $i:\R^q\to\R^n$ the embedding given by
$i(t_1,\dotsc,t_q)=(t_1,\dotsc,t_q,0,\dotsc,0)$. Fix $\epsilon>0$, such
that the critical points $x$ of $\omega$ admit disjoint admissible charts
(in which (a) and (b) are satisfied)
with $\epsilon_x>\epsilon$. For $c>0$ choose a smooth complex valued function
$a_c(t_1,\dotsc,t_q)$ with support in the disc of radius $\epsilon$ and
satisfying
$$
\int_{\R^q}e^{-sc(t^2_1+\cdots+t^2_q)}
a_c(t_1,\dotsc,t_q)dt^1\wedge\cdots\wedge dt^q=1,
\tag 5.3
$$
and a smooth function $\beta:\R_+\times\R_+\to\R_+$, so that
$\beta(\cdot,\lambda)$ has support equal to  $[0,\lambda]$ and
satisfies $\beta(t,\lambda)=1$ for $0\leq t\leq\lambda/2$. Denote by
$a_{c,\lambda}:\R^n\to\C$ the function defined by
$$
a_{c,\lambda}(t_1,\dotsc,t_n)=
\beta\Big(\sqrt{t^2_{q+1}+\cdots+t^2_n},\lambda\Big)a_c(t_1,\dotsc,t_q)
$$
and by $\alpha_{c,\lambda}\in\Omega^q(\R^n;\C)$ the smooth form given by
$$
\alpha_{c,\lambda}
=a_{c,\lambda}(t_1,\dotsc,t_n)dt^1\wedge\cdots\wedge dt^q.
$$
Since the support of $\alpha_{c,\lambda}$ is contained in $B(\varepsilon,0)$,
we can, for every $x\in\Cr_q(\omega)$, define
$\alpha^x_\lambda\in\Omega^q(M;\C)$
by $(\theta_x^{-1})^*\alpha_{c_x,\lambda}$ on $U_x$ and extend it by zero.

For every $x,x'\in\Cr_q(\omega)$, we consider the
function $a^{x,x'}_\lambda(t_1,\dotsc,t_q)$, defined by
$(i_{x'}^-)^*(\alpha^x_\lambda)
=a^{x,x'}_\lambda(t_1,\dotsc,t_q)dt^1\wedge\cdots\wedge dt^q$
and observe that it has the following properties:
\roster
\item
$a^{x,x}_\lambda(t_1,\dotsc,t_q)=a_{c_x}(t_1,\dotsc,t_q)$ for all $\lambda>0$
and all $t_1^2+\cdots+t_q^2\leq\epsilon^2$.
\item
For $\lambda\leq\lambda'$ one has
$$
\bigl|a^{x,x'}_\lambda(t_1,\dotsc,t_q)\bigr|
\leq\bigl|a^{x,x'}_{\lambda'}(t_1,\dotsc,t_q)\bigr|
$$
and $\supp(a^{x,x'}_\lambda)\subseteq\supp(a^{x,x'}_{\lambda'})$.
\item
If $x\neq x'$ then for any compact $K\subset\R^q$ there exists
$\lambda$ small enough so that
$\supp(a^{x,x'}_\lambda)\cap K=\emptyset$.
\item
If $x=x'$ then for every compact $K\subseteq\R^q\setminus B(\epsilon,0)$ there
exists $\lambda$ small enough so that
$\supp(a^{x,x}_\lambda)\cap K=\emptyset$.
\endroster
If $x=x'$ \therosteritem1, \therosteritem2, \therosteritem4 and \thetag{5.3}
imply
$$
\lim_{\lambda\to 0}\int_{W_{\tilde x}^-}e^{sh^{\tilde x}}
\pi^*(\alpha_\lambda^x)
=1+\lim_{\lambda\to 0}\int_{\R^q\setminus B(\epsilon,0)}
e^{sh^{\tilde x}}a^{x,x}_\lambda dt^1\wedge\cdots\wedge t^q=1,
$$
where we also used the fact, that the integrals converge and applied the
dominant convergence theorem. If $x\neq x'$ the
same argument but using now \therosteritem3 instead of \therosteritem4, yields
$$
\lim_{\lambda\to 0}\int_{W_{\tilde x'}^-}e^{sh^{\tilde x'}}
\pi^*(\alpha_\lambda^x)
=\lim_{\lambda\to 0}\int_{\R^q}
e^{sh^{\tilde x'}}a^{x,x'}_\lambda dt^1\wedge\cdots\wedge t^q
=0.
$$

In order to prove part~\therosteritem2 of Proposition~4 note first
that we have
$d(e^{sh^{\tilde x}}\pi^*\alpha)=e^{sh^{\tilde x}}\pi^*(d^q_s(\alpha))$.
So we have
$$
\int_{\hat W_{\tilde x}^-}\hat i_{\tilde x}^*
\bigl(e^{sh^{\tilde x}}\pi^*(d^q_s(\alpha))\bigr)
=
\int_{\hat W_{\tilde x}^-}d\,\hat i_{\tilde x}^*
(e^{sh^{\tilde x}}\pi^*\alpha)
=
\int_{(\hat W_{\tilde x}^-)_1}\hat i_{\tilde x}^*
(e^{sh^{\tilde x}}\pi^*\alpha).
\tag 5.4
$$
To check the second equality in \thetag{5.4} we proceed as follows. 
Consider a smooth function $\beta:\R\to[0,1]$ which satisfies 
$\beta(t)=1$, if $t\leq 0$, $\beta(t)=0$, if $t\geq 1$ and 
$-2\leq\beta'(t)\leq 0$. For any a positive integer $N$,
denote by $\rho_N:[0,\infty)\to[0,1]$ the function $\rho_N(t)=\beta(t-N)$.
Define the smooth function $\chi_N:\hat W_{\tilde x}^-\to[0,1]$ by
$\chi_N:=\rho_N\circ h^{\tilde x}\circ\hat i_{\tilde x}$. Clearly 
$\chi_N$ has compact support contained in $(h^{\tilde x})^{-1}([0,N+1])$, 
since $h^{\tilde x}\circ\hat i_{\tilde x}$ is proper, cf 
Theorem~1\therosteritem2.

Observe that
$$
\multline 
\int_{\hat W_{\tilde x}^-}d\,\hat i_{\tilde x}^*
(e^{sh^{\tilde x}}\pi^*\alpha)
=\\=\lim_{N\to\infty}
\int_{\hat W_{\tilde x}^-}d\bigl(\chi_N\hat i_{\tilde x}^*
(e^{sh^{\tilde x}}\pi^*\alpha)\bigr)
-\lim_{N\to\infty}
\int_{\hat W_{\tilde x}^-}\hat i_{\tilde x}^*\bigl((\rho'_N\circ h^{\tilde x})
e^{sh^{\tilde x}}\pi^*(\omega\wedge\alpha)\bigr)
\endmultline
\tag 5.5
$$
Note that 
$$
\int_{\hat W_{\tilde x}^-}\bigl|\hat i^*_{\tilde x}
\bigl((\rho'_N\circ h^{\tilde x})e^{sh^{\tilde x}}\pi^*(\omega\wedge\alpha)
\bigr)\bigr|
\leq2\int_{(h^{\tilde x}\circ\hat i_{\tilde x})^{-1}(N,N+1)}
\bigl|\hat i_{\tilde x}^*\bigl(e^{sh^{\tilde x}}
\pi^*(\omega\wedge\alpha)\bigr)\bigr|.
$$
Then, in view of the absolute convergence of 
$\int_{\hat W_{\tilde x}^-}\hat i_{\tilde x}^*
\bigl(e^{sh^{\tilde x}}\pi^*(\omega\wedge\alpha)\bigr)$, one concludes that 
the second limit in the right side of \thetag{5.5} is zero,
and therefore by Stoke's theorem we derive the second equality of 
\thetag{5.4}. 

The left hand side of \thetag{5.4} is
$$
\int_{\hat W_{\tilde x}^-}\hat i_{\tilde x}^*
\bigl(e^{sh^{\tilde x}}\pi^*(d^q_s(\alpha))\bigr)
=\Int_s(d^q_s(\alpha))(\tilde x).
$$
To compute the right side let $0>a_1>a_2>\cdots>a_k>\cdots$ be a sequence of 
regular values for $h^{\tilde x}$ restricted to $W_{\tilde x}$
tending to $-\infty$ and denote by $\hat W_{\tilde x}^-(n)$ the subset 
$(\hat W_{\tilde x}^-)_1\cap(h^{\tilde x})^{-1}([0,a_n])
=\bigcup_{\{\tilde y\in\Cal X_q\mid
\text{$\tilde y$ is $\tilde x$-relevant and $h^{\tilde x}(\tilde y)>a_n$}\}}
\Cal T(\tilde x,\tilde y)\times\hat W^-_{\tilde y}$.

Using the description of the boundary of $\hat W^-_{\tilde x}$ in
Theorem~1\therosteritem2, the  
convergence of the integrals $\int_{\Cal T(\tilde x,\tilde y)\times
\hat W^-_{\tilde y}}\hat i_{\tilde x}^*(e^{sh^{\tilde x}}\pi^*\alpha)$, 
(assured by Proposition~4\therosteritem1 and the finiteness of the set 
$\Cal T(\tilde x,\tilde y)$) and the dominant convergence theorem,
the right hand side of \thetag{5.4} gives
$$
\split
\int_{(\hat W_{\tilde x}^-)_1}\hat i_{\tilde x}^*
(e^{sh^{\tilde x}}\pi^*\alpha)
&= \lim_{n\to\infty}
\int_{\hat W_{\tilde x}^-(n)}\hat i_{\tilde x}^*
(e^{sh^{\tilde x}}\pi^*\alpha)\\&=
\sum_{
\Sb\tilde y\in\Cal X_q
\\
\text{$\tilde y$ is $\tilde x$-relevant}\endSb}
\int_{\Cal T(\tilde x,\tilde y)\times\hat W^-_{\tilde y}}
\hat i_{\tilde x}^*(e^{sh^{\tilde x}}\pi^*\alpha)
\\&=\sum_{\tilde y\in\Cal X_q}I_{q+1}(\tilde x,\tilde y)
\int_{\hat W^-_{\tilde y}}\hat i_{\tilde y}^*(e^{sh^{\tilde x}}\pi^*\alpha)
\\&=\sum_{\tilde y\in\Cal X_q}I_{q+1}(\tilde x,\tilde y)
\int_{W^-_{\tilde y}}e^{sh^{\tilde x}}\pi^*\alpha
\\&=\sum_{\tilde y\in\Cal X_q}I_{q+1}(\tilde x,\tilde y)
e^{-sH(\tilde x,\tilde y)}\int_{W_{\tilde y}^-}e^{sh^{\tilde y}}\pi^*\alpha
\\&=\sum_{\tilde y\in\Cal X_q}I_{q+1}(\tilde x,\tilde y)
e^{-sH(\tilde x,\tilde y)}\Int_s(\alpha)(\tilde y),
\endsplit
$$
where we used $h^{\tilde x}=h^{\tilde y}-H(\tilde x,\tilde y)$ for the
fifth equality.
\qed
\enddemo

We close the section with the

\demo{Proof of Theorem 2}
Part \therosteritem1 follows immediately from the fact that
$\gamma \Cal M(\tilde x,\tilde y)=\Cal M(\gamma\tilde x,\gamma\tilde y)$.
To check \therosteritem2 observe that in view of Theorem~1\therosteritem1
$\Cal B(\tilde x,\tilde y)$ is a compact oriented smooth manifold with
corners of dimension one hence a disjoint union of oriented closed intervals
and circles. It is not hard to see that the left side of \thetag{1.5}
is nothing but the algebraic cardinality of the boundary of
$\Cal B(\tilde x,\tilde y)$, which has to be zero.

To check \therosteritem3 let $\tilde y\in\Cal X_q$ and choose
$\alpha_{\tilde y}\in\Omega^q(M;\C)$, so
that $\Int_s(\alpha_{\tilde y})=\delta_{\pi(\tilde y)}$. This is possible in
view of the surjectivity stated in Proposition~4\therosteritem1. By
applying Proposition~4\therosteritem2 to the form $\alpha_{\tilde y}$ we get
for every $\tilde x\in\Cal X_{q+1}$
$$
\Int_s\bigl(d^q_s(\alpha_{\tilde y})\bigr)(\tilde x)
=\sum_{\gamma\in\Gamma}I_{q+1}(\tilde x,\gamma\tilde y)
e^{-sH(\tilde x,\gamma\tilde y)}
=e^{-sH(\tilde x,\tilde y)}
\sum_{\gamma\in\Gamma}I_{q+1}(\gamma\tilde x,\tilde y)e^{-s[\omega](\gamma)}.
$$
By Proposition~4\therosteritem1 the left hand side converges, and hence so
does \thetag{1.6}.
\qed
\enddemo

\head 6. Sketch of the Proof of Theorems 3 and 4\endhead

First observe that the Witten Laplacians $\Delta^q_t$ are zero order
perturbation of the Laplace Beltramy operator $\Delta^q=\Delta^q_0$.
Precisely, cf \cite{HeSj84},
$$
\Delta^q_t=\Delta^q+t\bigl(L_{\grad_g\omega}+L_{\grad_g\omega}^\sharp\bigr)
+t^2\|\omega\|^2\id,
\tag 6.1
$$
where $\|\omega\|=\llangle\omega,\omega\rrangle$,
$L_{\grad_g\omega}$ denotes the Lie derivative with respect to the
vector field $\grad_g\omega$ and
$L_{\grad_g\omega}^\sharp:\Omega^q(M)\to\Omega^q(M)$ its formal adjoint
$$
L_{\grad_g\omega}^\sharp\alpha=
(-1)^{nq+q+1}*L_{\grad_g\omega}(*\alpha)
=d^\sharp(\omega\wedge\alpha)+\omega\wedge d^\sharp\alpha.
$$
Despite the fact that $L_{\grad_g\omega}$ is an order
one differential operator the operator
$L_{\grad_g\omega}+L_{\grad_g\omega}^\sharp$ has order zero.

In the neighborhood of a critical point $y$ and with respect to a chart
$(\theta_y,\epsilon_y)$ which satisfies \thetag{a} and \thetag{b}
the Witten Laplacian $\Delta^q_t$ (denoted in this case $\Delta^q_{k,t}$ to
emphasize the dependence on the index $k$)
can be written down as
$$
\Delta^q_{k,t}=\Delta^q+2c_ytM_{q,k}+4c_y^2t^2(x_1^2+\cdots+x_n^2)\id,
$$
with
$$
\Delta^q\Big(\sum_Ia_I(x_1,x_2,\dotsc,x_n)dx^I\Big)
=-\sum_I\Big(\sum_{i=1}^n\frac{\partial^2 a_I}{\partial x_i^2}
(x_1,x_2,\dotsc,x_n)\Big)dx^I,
$$
and $M_{q,k}$ is the linear operator determined by
$$
M_{q,k}\Big(\sum_Ia_I(x_1,x_2,\dotsc,x_n)dx^I\Big)=
\sum_I\epsilon_I^ka_I(x_1,x_2,\dotsc,x_n)dx^I.
$$
Here $I=(i_1,i_2,\dotsc,i_q)$, $1\leq i_1<i_2\cdots<i_q\leq n$,
$dx^I=dx^{i_1}\wedge\cdots\wedge dx^{i_q}$ and
$$
\epsilon_I^k=-n+2k-2\big|I\cap\{1,\dotsc,k\}\big|
+2\big|I\cap\{k+1,\dotsc,n\}\big|,
$$
where $|A|$ denotes the cardinality of the set $A$.
Note that $\epsilon^k_I\geq -n$, and equals $-n$ iff $q=k$ and
$I=(1,\dotsc,q)$, cf \cite{BFKM96}, page~804.

The proof of Theorem~3 is is based on a mini-max
criterion for detecting a gap
in the spectrum of a positive selfadjoint operator in a Hilbert space
$H$, cf Lemma~4 below, and some basic estimates for the harmonic oscillator
collected in Lemma~5 and 6 below.

\proclaim{Lemma 4}
Let $A:H\to H$ be a densely defined (not necessary bounded)
self adjoint positive operator in a Hilbert space
$\bigl(H,\langle\cdot,\cdot,\rangle\bigr)$ and $a,b$ two real numbers so
that $0<a<b<\infty$. Suppose that there exist two closed subspaces $H_1$ and
$H_2$ of $H$ with $H_1\cap H_2=0$ and $H_1+H_2=H$, such that
\roster
\item
$\langle Ax_1,x_1\rangle\leq a\|x_1\|^2$ for any $x_1\in H_1$, and
\item
$\langle Ax_2,x_2\rangle\geq b\|x_2\|^2$ for any $x_2\in H_2$.
\endroster
Then $\spect(A)\cap(a,b)=\emptyset$. Moreover if $H_1$ is finite dimensional
then $\dim H_1$ equals the number of eigenvalues of $A$ which are
smaller than $a$, counted with multiplicity.
\endproclaim

The proof of this lemma is elementary, cf Lemma~1.2 in \cite{BFK98} or
the proof of Proposition~5.2 in \cite{BFKM96}, pages~806--807.

Let $\Cal S^q(\R^n)$ denote the space
of smooth $q$-forms
$\omega=\sum_Ia_I(x_1,x_2,\dotsc,x_n)dx^I$
with $a_I(x_1,x_2,\dotsc,x_n)$
rapidly decaying functions. The operator $\Delta^q_{k,t}$ acting on
$\Cal S^q(\R^n)$ is globally elliptic (in the sense of \cite{Sh87}
or \cite{H85}), selfadjoint and positive. This operator is the harmonic
oscillator in $n$ variables acting on $q$-forms and its properties can
be derived from the harmonic oscillator in one variable
$-\frac{d^2}{dx^2}+a+bx^2$ acting on functions. In particular the following
result holds.

\proclaim{Lemma 5}
Let $t>0$. Then:
\roster
\item
$\Delta^q_{k,t}$, regarded as an unbounded densely defined operator on
the $L^2$-com\-ple\-tion of $\Cal S^q(\R^n)$, is selfadjoint, positive and its
spectrum is contained in $4c_yt\N_0$, \ie positive integer multiples of
$4c_yt$.
\item
$\ker(\Delta^q_{k,t})=0$ if $k\neq q$, and $\dim\ker(\Delta^q_{q,t})=1$.
\item
Denote $|x|^2:=\sum_ix_i^2$. Then
$$
\omega_{q,t}=\bigl(\tfrac{2c_yt}\pi\bigr)^{n/4}
e^{-c_yt|x|^2}dx^1\wedge\cdots\wedge dx^q
$$
is the generator of $\ker(\Delta^q_{q,t})$ with $L^2$-norm $1$.
\endroster
\endproclaim

For a proof consult \cite{BFKM96}, page~806 (step~1 in the proof of
Proposition~5.2).

For $\eta>0$ choose a smooth function $\gamma_\eta:\R\to\R$, which satisfies
$$
\gamma_\eta(u)=
\cases
1&\text{if $u\leq\eta/2$, and}
\\
0&\text{if $u\geq\eta$.}
\endcases
$$
Introduce
$\tilde\omega_{q,t}^\eta\in\Omega^q_c(\R^n)$ defined by
$$
\tilde\omega_{q,t}^\eta(x):=\beta_{q,t}^{-1}\gamma_\eta(|x|)\omega_{q,t}(x),
$$
where
$$
\beta_{q,t}=\bigl(\tfrac{2c_yt}\pi\bigr)^{n/4}
\Big(\int_{\R^n}\gamma_\eta^2(|x|)e^{-2c_yt|x|^2}
dx^1\cdots dx^n\Big)^{1/2}.
$$
The smooth form $\tilde\omega_{q,t}^\eta$ has its support in the ball
$D_\eta(0)$, agrees with $\omega_{q,t}$ on the ball $D_{\eta/2}(0)$
and satisfies
$$
\bigl\langle\tilde\omega_{q,t}^\eta,\tilde\omega_{q,t}^\eta\bigr\rangle=1
\tag 6.2
$$
with respect to the scalar product $\langle\cdot,\cdot\rangle$ on
$\Cal S^q(\R^n)$, induced by the Euclidean metric. The following lemma
can be
obtained
by elementary calculations in coordinates in view of the
explicit
formula of
$\Delta^q_{k,t}$, cf~\cite{BFKM96}, Appendix~2.

\proclaim{Lemma 6}
For a fixed $r\in\N_0$ there
exist positive constants $C$, $C'$, $C''$,
$T_0$ and $\epsilon_0$, so that $t\geq T_0$ and
$\epsilon\leq\epsilon_0$ imply
\roster
\item
$$
\Big|\frac{\partial^{|I|}}{\partial x_I}\Delta^q_{q,t}(
\tilde\omega_{q,t}^\epsilon)(x)\Big|\leq Ce^{-C't},
$$
for any
$x\in\R^n$ and multi-index $I=(i_1,\dotsc,i_n)$ with
$|I|=i_1+\cdots+i_n\leq r$.
\item
$\langle\Delta^q_{k,t}\tilde\omega_{q,t}^{\epsilon},
\tilde\omega_{q,t}^{\epsilon}\rangle\geq 2t|q-k|$
\item
If $\alpha\perp\tilde\omega^\epsilon_{q,t}$ with respect to
the scalar product $\langle\cdot,\cdot\rangle$ then
$$
\langle\Delta^q_{q,t}\alpha,\alpha\rangle\geq C''t\|\alpha\|^2.
$$
\endroster
\endproclaim

For the proof of Theorems~3 and 4 we set the following notations.
We choose
$\epsilon>0$ so that for each $y\in\Cr(\omega)$ there exists an
admissible coordinate chart
$\theta_y:(U_y,y)\to(D_{\epsilon_y},0)$,
cf Definition~1, with
$\epsilon_y\geq2\epsilon$ and so that
$U_y\cap U_z=\emptyset$ for $y\neq z$.
Here we write $D_\rho$ for the disc $D_\rho(0)\subseteq\R^n$ of radius $\rho$
centered at $0$.

Choose once and for all such an admissible coordinate chart for each
$y\in\Cr_q(\omega)$.
Introduce the smooth forms
$\bar\omega_{y,t}\in\Omega^q(M)$
defined by
$$
\bar\omega_{y,t}|_{M\setminus\theta_y^{-1}(D_{2\epsilon})}:=0,
\quad\text{and}\quad
\bar\omega_{y,t}|_{\theta_y^{-1}(D_{2\epsilon})}:=
\theta_y^*(\tilde\omega^\epsilon_{q,t}).
$$
For any given $t>0$ the forms $\bar\omega_{y,t}\in\Omega^q(M)$,
$y\in\Cr_q(\omega)$, are orthonormal. Indeed, if $y,z\in\Cr_q(\omega)$,
$y\neq z$ then $\bar\omega_{y,t}$ and $\bar\omega_{z,t}$
have disjoint support, hence are orthogonal. Because the support
of $\bar\omega_{y,t}$ is contained in an admissible chart we have
$\langle\bar\omega_{y,t},\bar\omega_{y,t}\rangle=1$ by \thetag{6.2}.

For $t\geq T_0$, with $T_0$ given by Lemma~6, we introduce
the linear map
$$
J^q_t:\maps(\Cr_q(\omega),\R)\to\Omega^q(M),\quad
J^q_t(\delta_y):=\bar\omega_{y,t},
$$
where
$\delta_y\in\maps(\Cr(\omega),\R)$ is given by $\delta_y(z)=\delta_{y,z}$ for
$y,z\in\Cr(\omega)$.
$J^q_t$ is an isometry, for we have equipped
$\maps(\Cr(\omega),\R)$ with the scalar product which makes the base
$\delta_y$ orthonormal, thus in particular injective.

\demo{Proof of Theorem 3 (sketch)}
Take $H$ to be the
$L^2$-completion of
$\Omega^q(M)$
with respect to the
scalar product $\langle\cdot,\cdot\rangle$,
$H_1:=J^q_t\bigl(\maps(\Cr_q(\omega),\R)\bigr)$ and
$H_2=H_1^\perp$.

Let $T_0$, $C$, $C'$ and $C''$ be given by Lemma~6 and define
$$
C_1:=\inf_{z\in M'}\bigl\|\grad_g\omega(z)\bigr\|,
$$
where
$M'=M\setminus\bigcup_{y\in\Cr_q(\omega)}\theta_y^{-1}(D_\epsilon)$,
and
$$
C_2=\sup_{x\in M}\bigl\|\bigl(L_{\grad_g\omega}
+L_{\grad_g\omega}^\sharp\bigr)(z)\bigr\|.
$$
Here $\bigl\|\grad_g\omega(z)\bigr\|$ \resp $\bigl\|\bigl(L_{\grad_g\omega}+
L_{\grad_g\omega}^\sharp\bigr)(z)\bigr\|$ denotes
the norm of
the vector
$\grad_g\omega(z)\in T_zM$ \resp of the linear map
$$
\bigl(L_{\grad_g\omega}+L_{\grad_g\omega}^\sharp\bigr)(z):
\Lambda^q(T_z^*M)\to\Lambda^q(T_z^*M)
$$
with respect to the scalar product induced in $T_zM$ and
$\Lambda^q(T_z^*M)$ by $g(z)$. Recall that if $X$ is a vector field then
$L_X+L_X^\sharp$ is a zero order differential operator, hence an
endomorphism of the bundle $\Lambda^q(T^*M)\to M$.

We can use the constants $T_0$, $C$, $C'$, $C''$, $C_1$ and $C_2$ to
construct $C_3$ and $\epsilon_1$ so that for $t\geq T_0$
and $\epsilon\leq\epsilon_1$, we have
$\langle\Delta_q(t)\alpha,\alpha\rangle\geq
C_3t\langle\alpha,\alpha\rangle$ for any
$\alpha\in H_2$, cf \cite{BFKM96}, pages~808--810.

Now one can apply Lemma~4 whose hypotheses are satisfied for $a=Ce^{-C't}$,
$b=C3t$ and $t\geq T_0$. This finishes the proof of Theorem~3.
\qed
\enddemo

Let $Q^q_t$, $t\geq T_0$ denote the orthogonal projection in $H$ onto
$\Omega^q_\tsm(M)$, the span of the eigenvectors corresponding the
eigenvalues smaller than $1$. In view of the ellipticity of
$\Delta^q_t$ all these eigenvectors are smooth $q$-forms. An additional
important estimate is given by the following

\proclaim{Lemma 7}
For $r\in\N_0$ one can find
$\epsilon_0>0$ and $C_4,C_5$ so that
for $t\geq T_0$ as constructed above, and any $\epsilon\leq\epsilon_0$
one has, for any $f\in\maps(\Cr_q(M),\R)$ and any $0\leq p\leq r$,
$$
\bigl\|(Q^q_tJ^q_t-J^q_t)(f)\bigr\|_{C^p}\leq C_4e^{- C_5t}\|f\|,
$$
where $\|\cdot\|_{C^p}$ denotes the $C^p$-norm.
\endproclaim

The proof of Lemma~7 is contained in \cite{BZ92}, page~128 and
\cite{BFKM96}, page~811. Its proof requires \thetag{6.1},
Lemma~6 and general estimates coming from the ellipticity of $\Delta^q_t$.

\demo{Proof of Theorem 4 (sketch)}
Let $T_0$ be provided by Lemma~7. For $t\geq T_0$, let
$R^q_t$ be the isometry defined by
$$
R^q_t:=Q^q_tJ^q_t\bigl((Q^q_tJ^q_t)^\sharp Q^q_tJ^q_t\bigr)^{-1/2}:
\maps(\Cr_q(\omega),\R)\to\Omega^q_\tsm(M)
$$
and introduce
$E_{t,y}:=R^q_t(\delta_y)\in\Omega^q(M)$ for any $y\in\Cr_q(\omega)$.
Lemma~7 implies that there
exists $\epsilon>0$, $T_0$ and
$C$ so that for any $t\geq T_0$ and any $y\in\Cr_q(\omega)$ one has
$$
\sup_{z\in M\setminus\theta_y^{-1}(D_\epsilon)}
\|E_{t,y}(z)\|\leq Ce^{-\epsilon t}
\tag 6.3
$$
and
$$
\bigl\|E_{t,y}(z)-\bar\omega_{y,t}(z)\bigr\|
\leq C\tfrac1t,
\tag 6.4
$$
for any $z\in W_y^-\cap\theta_y^{-1}(D_\epsilon)$.
To check Theorem 4 it suffices to show that

$$
\bigl|\Int_t(E_{t,y})(z)-\delta_y(z)\bigr|\leq C''\tfrac1t,
$$
for some $C''>0$ and any $y,z\in\Cr_q(\omega)$. For $y=z$ this follows from
\thetag{6.3} and for $y\neq z$ it follows from \thetag{6.4}.
\qed
\enddemo

\enddocument